# FLUID MODEL FOR A NETWORK OPERATING UNDER A FAIR BANDWIDTH-SHARING POLICY

BY F. P. KELLY[1] AND R. J. WILLIAMS[1,2]

*University of Cambridge and University of California*

We consider a model of Internet congestion control that represents the randomly varying number of flows present in a network where bandwidth is shared fairly between document transfers. We study critical fluid models obtained as formal limits under law of large numbers scalings when the average load on at least one resource is equal to its capacity. We establish convergence to equilibria for fluid models and identify the invariant manifold. The form of the invariant manifold gives insight into the phenomenon of entrainment whereby congestion at some resources may prevent other resources from working at their full capacity.

**1. Introduction.** Roberts and Massoulié [19] have introduced and studied a flow-level model of Internet congestion control, that represents the randomly varying number of flows present in a network where bandwidth is dynamically shared between flows that correspond to continuous transfers of individual documents. This model assumes a "separation of time scales" such that the time scale of the flow dynamics (i.e., of document arrivals and departures) is much longer than the time scale of the packet level dynamics on which rate control schemes such as TCP converge to equilibrium.

Subsequent to the work of Roberts and Massoulié, assuming exponentially distributed document sizes, de Veciana, Lee and Konstantopoulos [9] and Bonald and Massoulié [2] studied the stability of the flow-level model operating under various bandwidth sharing policies, where a bandwidth sharing

Received March 2003; revised June 2003.
[1]Supported in part by the Operations, Information and Technology Program and the Center for Electronic Business and Commerce of the Graduate School of Business, Stanford University.
[2]Supported in part by NSF Grants DMS-00-71408 and DMS-03-05272, and a John Simon Guggenheim Fellowship.

policy corresponds to a generalization of the notion of a processor sharing discipline from a single resource to a network with several shared resources. Lyapunov functions constructed in [9] for weighted max–min fair and proportionally fair policies, and in [2] for weighted $\alpha$-fair policies [$\alpha \in (0, \infty)$] [17], imply positive recurrence of the Markov chain associated with the model when the average load on each resource is less than its capacity.

As a mechanism for performance analysis, we propose to use critical fluid models and related Brownian models to explore the behavior of flow-level models operating under weighted $\alpha$-fair bandwidth sharing policies in heavy traffic. We are particularly interested in manifestations of the phenomenon of entrainment, whereby congestion at some resources may prevent other resources from working at their full capacity. As a first step in this exploration, in this paper we consider critical fluid models, obtained as formal limits under law of large numbers scaling, from the flow-level models with exponentially distributed document sizes and operating under weighted $\alpha$-fair bandwidth sharing policies. The term critical refers to the fact that the nominal (or average) load on at least one resource is equal to its capacity, and for the other resources their nominal loads do not exceed their capacities, see (11) and (12). We identify the invariant states for the critical fluid models and we study the convergence to equilibria of critical fluid model solutions as time goes to infinity. Extrapolating from results for open multiclass queueing networks, we conjecture that such behavior is key to establishing heavy traffic diffusion approximations (also called Brownian models) for these flow-level models. We indicate the natural diffusion approximations suggested by our fluid model results.

There are several motivations for our work. One source of motivation lies in fixed point approximations of network performance for TCP networks (cf. [[7], [12], [20]]). These approximations require, as input, information on the joint distribution of the numbers of flows present on different routes, where dependencies between these numbers may be induced by the bandwidth sharing mechanism. Similarly, an understanding of such joint distributions seems important if the performance models for a single bottleneck described by Ben Fredj, Bonald, Proutiere, Regnie and Roberts [1] are to be generalized to a network. Another motivation is that the flow-level model typically involves the simultaneous use of several resources. With exponential document sizes, this model can be equated (in distribution) with a stochastic processing network (SPN) as introduced by Harrison [13, 14]. Open multiclass queueing networks are a special case of SPNs without simultaneous resource possession. For such networks operating under a head-of-the-line service discipline, it has been shown [5, 21] that suitable asymptotic behavior of critical fluid models implies a property called state space collapse, which validates the use of Brownian model approximations for these networks in heavy traffic. For more general SPNs, investigation of the behavior of critical fluid models, of a related notion of state space collapse, and of the implications for diffusion approximations, are in the early stages of development. The analysis in this paper can be viewed as a contribution to such an investigation for models involving simultaneous resource possession. Finally, although we restrict to exponential document sizes in this paper, we would like to relax that assumption in future work. Although this involves a significantly more elaborate stochastic model to keep track of residual document sizes (because of the processor sharing nature of the bandwidth sharing policy), knowing the results for exponential document sizes is likely to be useful for such work.

In order to state our results for the fluid model and conjectures for diffusion approximations, we need first to define the network structure, the weighted $\alpha$-fair bandwidth sharing policy and the stochastic model. This is done in Sections 2–4. The notion of a fluid model solution is defined in Section 5 and we state our main results there. The proofs of these results are given in Section 6. Appendix A develops some properties of the function that defines bandwidth allocations and Appendix B shows that our definition of a fluid model solution is reasonable in that fluid model solutions can be obtained as limit points of the stochastic model under fluid (or law of large numbers) scaling.

*Notation.* For each positive integer $d \geq 1$, $\mathbb{R}^d$ will denote $d$-dimensional Euclidean space and the positive orthant in this space will be denoted by $\mathbb{R}_+^d = \{x \in \mathbb{R}^d : x_i \geq 0 \text{ for } i = 1, \ldots, d\}$. The Euclidean norm of $x \in \mathbb{R}^d$ will be denoted by $\|x\|$. Inequalities between vectors in $\mathbb{R}^d$ will be interpreted componentwise, that is, for $x, y \in \mathbb{R}^d$, $x \leq y$ is equivalent to $x_i \leq y_i$ for $i = 1, \ldots, d$. Given a vector $x \in \mathbb{R}^d$, the $d \times d$ diagonal matrix with the entries of $x$ on its diagonal will be denoted by $\mathrm{diag}(x)$. For positive integers $d_1$ and $d_2$, the norm of a $d_1 \times d_2$ matrix $A$ will be given by

$$\|A\| = \left( \sum_{i=1}^{d_1} \sum_{j=1}^{d_2} A_{ij}^2 \right)^{1/2}.$$

The set of nonnegative integers will be denoted by $\mathbb{Z}_+$ and the set of points in $\mathbb{R}_+^d$ with all integer coordinates will be denoted by $\mathbb{Z}_+^d$. A sum over an empty set of indices will be taken to have a value of zero. The cardinality of a finite set $S$ will be denoted by $|S|$.

**2. Network structure.** We consider a network with finitely many *resources* labelled by $j \in \mathcal{J}$. A *route* $i$ is a nonempty subset of $\mathcal{J}$ (interpreted as the set of resources used by route $i$). We are given a set $\mathcal{I}$ of allowed routes. We assume that $\mathcal{J}$ and $\mathcal{I}$ are both nonempty and finite. Let $\mathbf{J} = |\mathcal{J}|$, the total number of resources, and $\mathbf{I} = |\mathcal{I}|$, the total number of routes. Let $A$ be the $\mathbf{J} \times \mathbf{I}$ matrix containing only zeros and ones, defined such that $A_{ji} = 1$ if resource $j$ is used by route $i$ and $A_{ji} = 0$ otherwise. Our assumption that each route $i$ identifies a *nonempty* subset of $\mathcal{J}$ implies that no column of $A$ is identically zero. We assume that $A$ has rank $\mathbf{J}$, so that it has full row rank. We further assume that *capacities* $(C_j : j \in \mathcal{J})$ are given and that these are all strictly positive and finite.

**3. Bandwidth sharing policy.** Bandwidth is allocated dynamically to the routes according to the following bandwidth sharing policy, which was first introduced by Mo and Walrand [17]. (To see how this fits into a stochastic model for the network dynamics, see Section 4.)

Given a fixed parameter $\alpha \in (0, \infty)$ and strictly positive weights $(\kappa_i : i \in \mathcal{I})$, if $N_i(t)$ denotes the (random) number of flows on route $i$ at time $t$ for each $i \in \mathcal{I}$ and $N(t) = (N_i(t) : i \in \mathcal{I})$, then the bandwidth allocated to route $i$ at time $t$ is given by $\Lambda_i(N(t))$ and this bandwidth is shared equally amongst all of the flows on route $i$. The function $\Lambda(\cdot) = (\Lambda_i(\cdot) : i \in \mathcal{I})$ is defined as follows (we define it on all of $\mathbb{R}_+^{\mathbf{I}}$ as we shall later apply it to fluid analogues of $N$). Let $\Lambda : \mathbb{R}_+^{\mathbf{I}} \to \mathbb{R}_+^{\mathbf{I}}$ be defined such that for each $n \in \mathbb{R}_+^{\mathbf{I}}$, $\Lambda_i(n) = 0$ for $i \in \mathcal{I}_0(n) \equiv \{l \in \mathcal{I} : n_l = 0\}$, and when $\mathcal{I}_+(n) \equiv \{l \in \mathcal{I} : n_l > 0\}$ is nonempty, $\Lambda^+(n) \equiv (\Lambda_i(n) : i \in \mathcal{I}_+(n))$ is the unique value of $\Lambda^+ = (\Lambda_i : i \in \mathcal{I}_+(n))$ that solves the optimization problem

(1)      maximize $G_n(\Lambda^+)$

(2)      subject to $\sum_{i \in \mathcal{I}_+(n)} A_{ji} \Lambda_i \leq C_j, \qquad j \in \mathcal{J},$

(3)      over    $\Lambda_i \geq 0, \qquad\qquad i \in \mathcal{I}_+(n),$

where for $n \in \mathbb{R}_+^{\mathbf{I}} \setminus \{0\}$ and $\Lambda^+ = (\Lambda_i : i \in \mathcal{I}_+(n)) \in \mathbb{R}_+^{|\mathcal{I}_+(n)|}$,

(4)    $G_n(\Lambda^+) = \begin{cases} \displaystyle\sum_{i \in \mathcal{I}_+(n)} \kappa_i n_i^\alpha \frac{\Lambda_i^{1-\alpha}}{1-\alpha}, & \text{if } \alpha \in (0, \infty) \setminus \{1\}, \\ \displaystyle\sum_{i \in \mathcal{I}_+(n)} \kappa_i n_i \log \Lambda_i, & \text{if } \alpha = 1, \end{cases}$

and the value of the right member above is taken to be $-\infty$ if $\alpha \in [1, \infty)$ and $\Lambda_i = 0$ for some $i \in \mathcal{I}_+(n)$. The resulting allocation is called a *weighted $\alpha$-fair allocation*.

Various properties of the mapping $\Lambda : \mathbb{R}_+^{\mathbf{I}} \to \mathbb{R}_+^{\mathbf{I}}$ are developed in Appendix A of this paper. In particular, for each $n \in \mathbb{R}_+^{\mathbf{I}}$:

(i) $\Lambda_i(n) > 0$ for each $i \in \mathcal{I}_+(n)$,
(ii) $\Lambda(rn) = \Lambda(n)$ for each $r > 0$,
(iii) $\Lambda_i(\cdot)$ is continuous at $n$ for each $i \in \mathcal{I}_+(n)$, and

(iv) there is $p \in \mathbb{R}_+^{\mathbf{J}}$ (not necessarily unique and depending on $n$) such that

$$\Lambda_i(n) = n_i \left( \frac{\kappa_i}{\sum_{j \in \mathcal{J}} p_j A_{ji}} \right)^{1/\alpha} \qquad \text{for each } i \in \mathcal{I}_+(n), \tag{5}$$

where

$$p_j \left( C_j - \sum_{i \in \mathcal{I}} A_{ji} \Lambda_i(n) \right) = 0 \qquad \text{for all } j \in \mathcal{J}. \tag{6}$$

The $(p_j : j \in \mathcal{J})$ are Lagrange multipliers for the optimization problem (1)–(3), one for each of the capacity constraints (2). [Note that for each $i \in \mathcal{I}_+(n)$, since $\Lambda_i(n) > 0$ [by (i)] and $n_i > 0$ (by definition), the fact that the representation (5) holds implies that $p$ is such that the denominator in the right member of (5) does not vanish.]

When $\kappa_i = 1, i \in \mathcal{I}$, the cases $\alpha \to 0$, $\alpha \to 1$ and $\alpha \to \infty$ correspond respectively to an allocation which achieves maximum throughput, is *proportionally fair* or is *max–min fair* [2, 17]. Weighted $\alpha$-fair allocations provide a tractable theoretical abstraction of decentralized packet-based congestion control algorithms such as TCP, the transmission control protocol of the Internet. Indeed, if $\alpha = 2$ and $\kappa_i$ is the reciprocal of the square of the round trip time on route $i$, then the formula (5) is a version of the *inverse square root law* familiar from studies of the throughput of TCP connections [11, 16, 18]. The relations (2) and (3), (5) and (6) and more refined versions of these relations, can be solved by iterative methods to give predictions of throughput, given the numbers of flows $N(t)$ present at time $t$ [7, 12, 20]. Given a distribution for $N(t)$, the overall network performance can be predicted. But a major difficulty with this approach is the choice of the distribution for $N(t)$. For example, if flows arrive on different routes as independent Poisson processes and if flows on a route remain in the system for independent and identically distributed holding periods, then the stationary distribution of the process $N$ is easy to describe: the components are independent, each with a Poisson distribution, whatever the distribution of holding periods. This model is indeed used in [12] and might be appropriate for real-time flows whose time in the system is unaffected by their allocated bandwidth. But for many flows, for example, document transfers, their length of time in the system is affected by their allocated bandwidth, and this may produce correlations between the components of $N(t)$ which need to be understood. Roberts and Massoulié [19] have begun the study of a stochastic model that captures this effect.

**4. Stochastic model.** An active flow on route $i$ corresponds to the continuous transmission of a document through the resources used by route $i$. Transmission is assumed to occur simultaneously through all resources on route $i$. The number of active flows on route $i$ at time $t$ is denoted by $N_i(t)$. The stochastic process $N = \{(N_1(t), \ldots, N_\mathbf{I}(t)), t \geq 0\}$ is assumed to be a Markov process with state space $\mathbb{Z}_+^\mathbf{I}$ and infinitesimal transition rates $q : \mathbb{Z}_+^\mathbf{I} \times \mathbb{Z}_+^\mathbf{I} \to \mathbb{R}$ given by

$$q(n, m) = \nu_i \qquad \text{if } m = n + e_i, \tag{7}$$

$$q(n, m) = \mu_i \Lambda_i(n) \qquad \text{if } m = n - e_i, \ n_i \geq 1, \tag{8}$$

$$q(n, m) = 0 \qquad \text{otherwise}, \tag{9}$$

for each $n, m \in \mathbb{R}_+^\mathbf{I}$, $i \in \mathcal{I}$, where, for each $i$, $\nu_i > 0$ and $\mu_i > 0$ are fixed constants, and $e_i$ is the **I**-dimensional unit vector whose $i$th component is 1 and whose other components are all zero.

This corresponds to a model where new flows arrive on route $i$ according to a Poisson process of rate $\nu_i$; for $i$ such that $N_i(t) \neq 0$, $\Lambda_i(N(t))/N_i(t)$ is the bandwidth allocated to each active flow on route $i$ at time $t$; and a flow on route $i$ transfers a document whose size is exponentially distributed with parameter $\mu_i$. This is the model of Roberts and Massoulié [19] with exponential document sizes.

From the results of de Veciana, Lee and Konstantopoulos [9] and Bonald and Massoulié [2], we know that the Markov chain $N$ is positive recurrent if

$$\sum_{i \in \mathcal{I}} A_{ji} \rho_i < C_j, \qquad j \in \mathcal{J}, \tag{10}$$

where $\rho_i = \nu_i/\mu_i$ for all $i \in \mathcal{I}$. These are natural constraints: $\rho_i$ is the average load produced by route $i$, and we can identify the ratio of the two sides of the inequality (10) as the *traffic intensity* at resource $j$. Indeed, condition (10) is necessary for positive recurrence of $N$. For a proof, suppose that $N$ is positive recurrent and fix $j \in \mathcal{J}$. The virtual waiting time $V_j(t)$ for resource $j$ at time $t$ is the amount of time, measured from time $t$ onwards, that it would take to complete the transfer of all of the documents that are being transmitted through resource $j$ at time $t$, assuming that external arrivals are turned off after time $t$, that is, no new documents are accepted for transmission after time $t$, and that all other resources are given infinite capacity, that is, $C_k = +\infty$ for all $k \neq j$, after time $t$. The virtual waiting time thus measures the time it would take for resource $j$ to become idle if there were no more arrivals after time $t$ and if resource $j$ could work at full capacity from time $t$. Suppose that the network starts empty. The positive recurrence of $N$ implies that the mean time for the virtual waiting time process $V_j$ to return to zero (after first moving away from zero) is finite. Consider another network with the same features as the original one, except that $C_k = +\infty$ for all $k \neq j$. Let $\tilde{V}_j$ denote the virtual waiting time process for resource $j$ in this network. When the same arrival and document size processes are used for the two networks, $\tilde{V}_j(t) \leq V_j(t)$ for all $t$. In particular, the mean time for $\tilde{V}_j$ to return to zero must be finite. Now, $\tilde{V}_j$ is equivalent in distribution to the virtual waiting time process for a multiclass single server queueing system operating under a work conserving service discipline. This system has one queue for each $i$ such that $A_{ji} = 1$. The queue associated with such an $i$ has an infinite capacity buffer, Poisson arrivals at rate $\nu_i$, i.i.d. exponential service times with a mean of $1/\mu_i$, and the server serves at a maximum rate of $C_j$. The virtual waiting time process for this queueing system is the same for all work conserving service disciplines and it is well known that the mean time for this process to return to zero is finite if and only if $\sum_{i \in \mathcal{I}} A_{ji} \rho_i < C_j$. Since $j$ was arbitrary, it follows that (10) must hold.

It is an open question whether, in the generalization of the above model to allow arbitrarily (rather than exponentially) distributed document sizes, the condition (10) is sufficient for stability.

**5. Main results.** Our aim in this paper is to begin to explore the behavior of the Markov chain $\{N(t), t \geq 0\}$ when

$$\text{(11)} \qquad \sum_{i \in \mathcal{I}} A_{ji} \rho_i \leq C_j, \qquad j \in \mathcal{J},$$

and some of the constraints are saturated, that is, some of the resources are in heavy traffic. Thus, we henceforth assume that (11) holds and that

$$\text{(12)} \qquad \mathcal{J}_* \equiv \left\{ j \in \mathcal{J} : \sum_{i \in \mathcal{I}} A_{ji} \rho_i = C_j \right\} \neq \varnothing.$$

Let $\mathbf{J}_* = |\mathcal{J}_*|$ and without loss of generality assume that the first $|\mathbf{J}_*|$ elements of $\mathcal{J}$ correspond to the set $\mathcal{J}_*$.

Here, we focus on understanding the behavior of fluid model solutions, which can be thought of as formal limits of the stochastic process $N$ under law of large numbers scaling. The following notions are used in the definition below. A function $f = (f_1, \ldots, f_\mathbf{I}) : [0, \infty) \to \mathbb{R}_+^\mathbf{I}$ is absolutely continuous if each of its components $f_i : [0, \infty) \to \mathbb{R}_+$, $i = 1, \ldots, \mathbf{I}$, is absolutely continuous. A *regular point* for an absolutely continuous function $f : [0, \infty) \to \mathbb{R}_+^\mathbf{I}$ is a value of $t \in [0, \infty)$ at which each component of $f$ is differentiable. [Since $f$ is absolutely continuous, almost every time $t \in [0, \infty)$ is a regular point for $f$.]

DEFINITION 5.1. A fluid model solution is an absolutely continuous function $n : [0, \infty) \to \mathbb{R}_+^\mathbf{I}$ such that at each regular point $t$ for $n(\cdot)$, we have for each $i \in \mathcal{I}$,

$$\text{(13)} \qquad \frac{d}{dt} n_i(t) = \begin{cases} \nu_i - \mu_i \Lambda_i(n(t)), & \text{if } n_i(t) > 0, \\ 0, & \text{if } n_i(t) = 0, \end{cases}$$

and for each $j \in \mathcal{J}$,

$$(14) \qquad \sum_{i \in \mathcal{I}_+(n(t))} A_{ji} \Lambda_i(n(t)) + \sum_{i \in \mathcal{I}_0(n(t))} A_{ji} \rho_i \leq C_j,$$

where $\mathcal{I}_+(n(t)) = \{i \in \mathcal{I} : n_i(t) > 0\}$ and $\mathcal{I}_0(n(t)) = \{i \in \mathcal{I} : n_i(t) = 0\}$.

Motivation for this definition is given in Appendix B through a fluid limit result. For the moment, we observe that if $n_i(t) > 0$, then the right-hand side of (13) is the infinitesimal drift of $N_i(t)$ when $N_i(t) > 0$ [cf. (7) and (8)]. On the other hand, if $n_i(t) = 0$ and $t$ is a regular point for $n(\cdot)$, then the derivative of $n_i(\cdot)$ at $t$ is forced to be zero since $n_i(s) \geq 0$ for all $s \geq 0$—to see this, consider the left- and right-hand derivatives of $n_i(\cdot)$ at $t$. This property may seem counterintuitive, however, this phenomenon is common in fluid models for queueing systems. It reflects the fact that a fluid model solution is obtained as a (formal) law of large numbers limit from the original stochastic model, and consequently a fluid model solution state $n \in \mathbb{R}_+^\mathbf{I}$, for which $n_i = 0$ can be the limit of rescaled states in the stochastic model where the $i$th component is at or near zero. The inequality (14) is derived from the fact that, in the stochastic model, the cumulative unused capacity for each resource is a nondecreasing process. As in the derivation of the differential equation (13), some care is needed here in treating routes $i$, for which $n_i(t) = 0$. One might paraphrase (14) as saying that the total fluid model bandwidth allocation for each resource cannot exceed its capacity, where the allocation to any route $i$ satisfying $n_i(t) = 0$ is $\rho_i$ at time $t$. For a more detailed justification, we refer the reader to Theorem B.1 and its proof.

Following Bramson [5], we now define an invariant manifold for fluid model solutions.

DEFINITION 5.2. A state $n_0 \in \mathbb{R}_+^\mathbf{I}$ is called invariant if there is a fluid model solution $n(\cdot)$ such that $n(t) = n_0$ for all $t \geq 0$. Let $\mathcal{M}_\alpha$ denote the set of all invariant states. We call $\mathcal{M}_\alpha$ the invariant manifold.

The following is a simple characterization of $\mathcal{M}_\alpha$.

LEMMA 5.1. *The set of invariant states, $\mathcal{M}_\alpha$, is given by*

$$(15) \qquad \{n \in \mathbb{R}_+^\mathbf{I} : \Lambda_i(n) = \rho_i \text{ for all } i \in \mathcal{I}_+(n)\}.$$

PROOF. Let $\mathcal{N}_\alpha$ denote the set in (15). Note that $\mathcal{M}_\alpha$ and $\mathcal{N}_\alpha$ are nonempty since they both contain the origin in $\mathbb{R}_+^\mathbf{I}$.

To show that $\mathcal{M}_\alpha \subset \mathcal{N}_\alpha$, suppose that $n_0 \in \mathcal{M}_\alpha$. If $n_0 = 0$, then it follows trivially that $n_0 \in \mathcal{N}_\alpha$. If $n_0 \neq 0$, then there is a fluid model solution $n(\cdot)$ satisfying $n(t) = n_0$ for all $t \geq 0$ and so it follows from (13) that

(16) $$\Lambda_i(n_0) = \rho_i \qquad \text{for all } i \in \mathcal{I}_+(n_0).$$

Conversely, to show that $\mathcal{N}_\alpha \subset \mathcal{M}_\alpha$, suppose that $n_0 \in \mathcal{N}_\alpha$. Then, by (11), $n(t) = n_0$ for all $t \geq 0$ satisfies (13) and (14) for all $t$ and all $i \in \mathcal{I}, j \in \mathcal{J}$, and so $n(\cdot)$ is a valid fluid model solution. Hence $n_0$ is in $\mathcal{M}_\alpha$. □

The following alternative characterization of the invariant states will also be used. It is proved in Section 6.

THEOREM 5.1. *A state $n \in \mathbb{R}_+^{\mathbf{I}}$ is an invariant state if and only if there is $q \in \mathbb{R}_+^{\mathbf{J}_*}$ such that*

(17) $$n_i = \rho_i \left( \frac{\sum_{j \in \mathcal{J}_*} q_j A_{ji}}{\kappa_i} \right)^{1/\alpha} \qquad \text{for all } i \in \mathcal{I}.$$

REMARK 5.1. In fact, as examination of the proof of the above theorem reveals, for an invariant state $n$, there is a one-to-one correspondence between the vectors $q$ appearing in the above representation of $n$ and the Lagrange multipliers $p$ appearing in the characterization of $\Lambda^+(n)$ given in Lemma A.4 [see also (5) and (6)]. This correspondence is obtained by taking the entries of $q$ to be given by the entries of $p$ with indices $j \in \mathcal{J}_*$ and noting that the other entries in $p$ are necessarily zero.

For each $n \in \mathbb{R}_+^{\mathbf{I}}$, we define the distance of $n$ from $\mathcal{M}_\alpha$ as

(18) $$d(\mathcal{M}_\alpha, n) = \inf\{\|v - n\| : v \in \mathcal{M}_\alpha\}.$$

The following theorem shows that, starting in any compact set, fluid model solutions converge uniformly towards the invariant manifold. This theorem is proved in Section 6.

THEOREM 5.2. *Fix $R \in (0, \infty)$ and $\varepsilon > 0$. There is a constant $T_{R,\varepsilon} < \infty$ such that for each fluid model solution $n(\cdot)$ satisfying $\|n(0)\| \leq R$ we have*

(19) $$d(\mathcal{M}_\alpha, n(t)) < \varepsilon \qquad \text{for all } t > T_{R,\varepsilon}.$$

In the course of proving Theorem 5.2, in Section 6, we prove the following (see Theorem 5.3) alternative characterization of invariant states. For this, define $w(n) = (w_j(n) : j \in \mathcal{J}_*)$ for $n \in \mathbb{R}_+^{\mathbf{I}}$ to be given by

(20) $$w_j(n) = \sum_{i \in \mathcal{I}} A_{ji} \frac{n_i}{\mu_i}, \qquad j \in \mathcal{J}_*.$$

We call $w(n)$ the *workload* associated with $n$. Let

$$(21) \qquad F(n) = \frac{1}{\alpha+1} \sum_{i \in \mathcal{I}} \nu_i \kappa_i \mu_i^{\alpha-1} \left(\frac{n_i}{\nu_i}\right)^{\alpha+1} \qquad \text{for all } n \in \mathbb{R}^{\mathbf{I}}_+.$$

This function $F$ was used in [2] as a Lyapunov function to show positive recurrence of $N$ under the conditions (10). An intuitive interpretation of the function $F$ is as follows. If the the number of flows on each route is fixed and given by the components of $n \in \mathbb{R}^{\mathbf{I}}_+$, then by Little's law $n_i/\nu_i$ is the mean time that a flow on route $i$ spends in the system and the time a flow on route $i$ spends in the system is exponentially distributed with mean $n_i/\nu_i$. The $(\alpha+1)$st moment of this random variable is $\Gamma(\alpha+2)(n_i/\nu_i)^{\alpha+1}$, where $\Gamma(\cdot)$ is the usual Gamma function. Thus, given $n$, $F(n)$ can be interpreted as a weighted sum over the routes, where for route $i$, the summand is the weight $\nu_i \kappa_i \mu_i^{\alpha-1}/(\alpha+1)\Gamma(\alpha+2)$ times the $(\alpha+1)$st moment of the amount of time spent in the system by a flow on that route.

For $w \in \mathbb{R}^{\mathbf{J}_*}_+$, define $\Delta(w)$ to be the unique value of $n \in \mathbb{R}^{\mathbf{I}}_+$ that solves the following optimization problem:

$$(22) \qquad \begin{aligned} & \text{minimize } F(n) \\ & \text{subject to } \sum_{i \in \mathcal{I}} A_{ji} \frac{n_i}{\mu_i} \geq w_j, \qquad j \in \mathcal{J}_*, \\ & \text{over } n_i \geq 0, \qquad\qquad\qquad\quad i \in \mathcal{I}. \end{aligned}$$

REMARK 5.2. Since $A$ has full row rank and its only entries are zeros and ones, for each $w \in \mathbb{R}^{\mathbf{J}_*}_+$, the feasible set of the optimization problem (22) is nonempty, and then since $F$ is nonnegative on $\mathbb{R}^{\mathbf{I}}_+$ and $F(n) \to \infty$ as $\|n\| \to \infty$, (22) has an optimal solution. By the strict convexity of $F$, this solution is unique.

THEOREM 5.3. *A vector $n \in \mathbb{R}^{\mathbf{I}}_+$ is an invariant state if and only if $n = \Delta(w(n))$.*

The map $\Delta : \mathbb{R}^{\mathbf{J}_*}_+ \to \mathbb{R}^{\mathbf{I}}_+$ plays an analogous role for the flow-level model of [19] to the lifting maps occuring in Bramson's work [3, 4] on the asymptotic behavior of fluid models associated with multiclass queueing networks operating under certain head-of-the-line service disciplines. It is natural to conjecture that one might prove a state space collapse theorem for the flow-level model in an analogous manner to that in [5], and extend the diffusion approximation results developed for multiclass queueing networks in [21], to prove a diffusion approximation for the flow-level model. This suggests that, under suitable rescaling and initial conditions, a diffusion approximation for the $\mathbf{J}_*$-dimensional workload process $W = \{W(t) : t \geq 0\}$ defined by

$$W_j(t) = \sum_{i \in \mathcal{I}} A_{ji} \frac{N_i(t)}{\mu_i}, \qquad j \in \mathcal{J}_*, \tag{23}$$

is likely to be a reflecting Brownian motion $\tilde{W}$ living in the *workload cone*

$$\mathcal{W}_\alpha = A_* M^{-1} \mathcal{M}_\alpha, \tag{24}$$

where $M = \text{diag}(\mu)$, $A_*$ is the $\mathbf{J}_* \times \mathbf{I}$ matrix obtained from $A$ by eliminating those rows of $A$ that are not indexed by elements of $\mathcal{J}_*$, and

$$\mathcal{M}_\alpha = \left\{ n \in \mathbb{R}_+^{\mathbf{I}} : n_i = \rho_i \left( \frac{\sum_{j \in \mathcal{J}_*} q_j A_{ji}}{\kappa_i} \right)^{1/\alpha}, \right.$$

$$\left. i \in \mathcal{I}, \text{ for some } q \in \mathbb{R}_+^{\mathbf{J}_*} \right\}. \tag{25}$$

Here the direction of reflection on the boundary surface corresponding to $q_j = 0$ is the unit vector pointing in the direction of the positive $j$th coordinate axis. Furthermore, state space collapse should yield an approximation $\tilde{N}$ for $N$, under diffusion scaling, where $\tilde{N} = \Delta(\tilde{W})$. This conjecture will be pursued in a subsequent work.

To illustrate the conjecture, we consider the following simple example. Suppose that $\mathcal{J} = \{1, 2\}$ and $\mathcal{I} = \{\{1\}, \{2\}, \{1, 2\}\}$, corresponding to a linear network with two resources and three routes. Let $\alpha \in (0, \infty)$, $\kappa_i = \mu_i = 1$, for $i = 1, 2, 3$, $C_j = 1$ for $j = 1, 2$, and $\rho_1 + \rho_3 = \rho_2 + \rho_3 = 1$. Then the state space for the diffusion $\tilde{W}$ is the cone

$$\mathcal{W}_\alpha = \{(w_1, w_2) : w_1 = \rho_1 q_1^{1/\alpha} + \rho_3 (q_1 + q_2)^{1/\alpha},$$

$$w_2 = \rho_2 q_2^{1/\alpha} + \rho_3 (q_1 + q_2)^{1/\alpha}, \text{ for some } q_1 \geq 0, \ q_2 \geq 0\},$$

which, for all $\alpha \in (0, \infty)$, is the same as the cone

$$\{(w_1, w_2) : w_1 \geq 0, \ w_1 \rho_3 \leq w_2 \leq w_1 \rho_3^{-1}\}$$

pictured in Figure 1. Reflection occurs in the horizontal direction (corresponding to resource 1 incurring idleness) on the bounding face $w_1 = w_2 \rho_3$. The interpretation of this is that although there is work for resource 1 within the system, congestion at resource 2 is preventing resource 1 from working at its full capacity. Similarly, vertical reflection (corresponding to resource 2 incurring idleness) on the bounding face $w_2 = w_1 \rho_3$ is interpreted to mean that congestion at resource 1 is preventing resource 2 from working at its full capacity. Although the workload cone is the same for all $\alpha \in (0, \infty)$ in this example, this will not be the case, in general, for higher-dimensional workloads.

**6. Proofs: characterization of invariant states and convergence to the invariant manifold.**

PROOF OF THEOREM 5.1. It follows from Lemma 5.1 and the characterization of $\Lambda^+(n)$ in terms of Lagrange multipliers given in Lemma A.4, that $n \in \mathcal{M}_\alpha$ if and only if there is $p \in \mathbb{R}_+^{\mathbf{J}}$ such that

$$(26) \qquad n_i = \rho_i \left( \frac{\sum_{j \in \mathcal{J}} p_j A_{ji}}{\kappa_i} \right)^{1/\alpha} \qquad \text{for all } i \in \mathcal{I}_+(n)$$

and for all $j \in \mathcal{J}$,

$$(27) \qquad p_j \left( C_j - \sum_{i \in \mathcal{I}_+(n)} A_{ji} \rho_i \right) = 0.$$

Note that for $j \in \mathcal{J} \setminus \mathcal{J}_*$, $\sum_{i \in \mathcal{I}} A_{ji} \rho_i < C_j$ and so (27) holds for such a $j$ if and only if $p_j = 0$. It follows that we can replace $\mathcal{J}$ by $\mathcal{J}_*$ and $\mathbf{J}$ by $\mathbf{J}_*$ in the above characterization of invariant states. The characterization given in the theorem can now be deduced as follows.

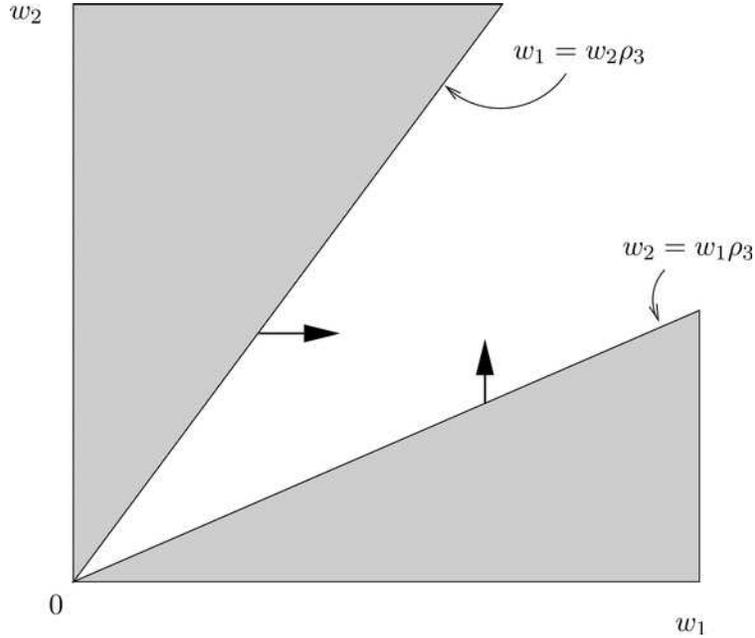

FIG. 1. The workload cone $\mathcal{W}_\alpha$ for a network with two resources, with workloads labelled $w_1, w_2$, and three routes, with traffic loads labelled $\rho_1, \rho_2, \rho_3$. Under the lifting map $\Delta$, points $(w_1, w_2)$ on the boundary $w_1 = w_2 \rho_3$ are mapped to points $(n_1, n_2, n_3)$, where $n_1 = 0$ (and the corresponding $q \in \mathbb{R}_+^2$ has $q_1 = 0$); similarly, points $(w_1, w_2)$ on the boundary $w_2 = w_1 \rho_3$ are mapped to points $(n_1, n_2, n_3)$, where $n_2 = 0$ (and the corresponding $q \in \mathbb{R}_+^2$ has $q_2 = 0$).

First, consider an invariant state $n$ and $i \in \mathcal{I} \setminus \mathcal{I}_+(n)$. Then $n_i = 0$ and for any $j \in \mathcal{J}_*$ such that $A_{ji} > 0$, we have

$$(28) \qquad \sum_{l \in \mathcal{I}_+(n)} A_{jl} \rho_l < \sum_{l \in \mathcal{I}} A_{jl} \rho_l = C_j,$$

and so by (27), we must have $p_j = 0$. On combining the above, we see that

$$(29) \qquad n_i = \rho_i \left( \frac{\sum_{j \in \mathcal{J}_*} p_j A_{ji}}{\kappa_i} \right)^{1/\alpha} \qquad \text{for all } i \in \mathcal{I}.$$

Thus, any invariant state has the form given in (17) with $q_j = p_j$ for $j \in \mathcal{J}_*$.

Conversely, suppose that $n$ is of the form given in (17) for some $q \in \mathbb{R}_+^{\mathbf{J}_*}$. Set $p_j = 0$ for $j \in \mathcal{J} \setminus \mathcal{J}_*$ and $p_j = q_j$ for $j \in \mathcal{J}_*$. Then, (26) holds immediately with $\mathcal{J}_*$ in place of $\mathcal{J}$, and so it suffices to show that $p$ and $n$ satisfy the complementarity condition (27) for each $j \in \mathcal{J}_*$. The only way that this can fail to hold is if there is $j \in \mathcal{J}_*$ and $i \in \mathcal{I}$ such that $q_j > 0$, $A_{ji} > 0$ and $n_i = 0$. But, by the representation (17), $n_i = 0$ implies that $q_j = 0$ for all $j \in \mathcal{J}_*$ satisfying $A_{ji} > 0$. Thus, (27) must hold for all $j \in \mathcal{J}_*$. □

We need some preliminary lemmas before we can prove the other characterization of invariant states given by Theorem 5.3. For this, recall the definition of the function $F$ from (21) in Section 5. For any fluid model solution $n(\cdot)$, $F(n(\cdot))$ is absolutely continuous and at each regular point $t$ for $n(\cdot)$,

$$(30) \qquad \frac{d}{dt} F(n(t)) = \sum_{i \in \mathcal{I}} \frac{\partial F}{\partial n_i} \frac{d}{dt} n_i(t)$$

$$(31) \qquad = \sum_{i \in \mathcal{I}_+(n(t))} \frac{\kappa_i}{\nu_i} \left( \frac{\mu_i}{\nu_i} \right)^{\alpha-1} (n_i(t))^{\alpha} (\nu_i - \mu_i \Lambda_i(n(t)))$$

$$(32) \qquad = \sum_{i \in \mathcal{I}_+(n(t))} \kappa_i \left( \frac{\mu_i n_i(t)}{\nu_i} \right)^{\alpha} \left( \frac{\nu_i}{\mu_i} - \Lambda_i(n(t)) \right)$$

$$(33) \qquad = K(n(t)),$$

where, for each $n \in \mathbb{R}_+^{\mathbf{I}}$,

$$(34) \qquad K(n) = \sum_{i \in \mathcal{I}_+(n)} \kappa_i \left( \frac{\mu_i n_i}{\nu_i} \right)^{\alpha} \left( \frac{\nu_i}{\mu_i} - \Lambda_i(n) \right).$$

Indeed, $K(0) = 0$ and for $n \in \mathbb{R}_+^{\mathbf{I}} \setminus \{0\}$,

$$(35) \qquad K(n) = \nabla G_n(\Lambda^{+,*}(n)) \cdot (\Lambda^{+,*}(n) - \Lambda^+(n)),$$

where $\Lambda^+(n) = (\Lambda_i(n) : i \in \mathcal{I}_+(n))$, $\Lambda^{+,*}(n) = (\rho_i : i \in \mathcal{I}_+(n))$, $\nabla G_n(\cdot)$ is the gradient of the function $G_n(\cdot)$ defined in (4).

LEMMA 6.1. *The function $K$ is continuous on $\mathbb{R}_+^{\mathbf{I}}$ and*

(36) $$K(n) \leq 0 \quad \text{for each } n \in \mathbb{R}_+^{\mathbf{I}},$$

*where the inequality is strict unless $n$ is an invariant state.*

PROOF. The continuity of $K$ follows from the definition (34), combined with Lemma A.3 and the fact that the term indexed by $i$ in the sum in (34) is small if $n_i$ is near zero, since $\Lambda(n)$ is bounded.

If $n = 0$, $K(n) = 0$ and 0 is an invariant state. Now suppose that $n \neq 0$. Then, by Lemma A.1, $\Lambda^+(n)$ solves the optimization problem (1)–(3) on $(0, \infty)^{|\mathcal{I}_+(n)|}$. Since $\Lambda^+ = \Lambda^{+,*}(n)$ is feasible for this problem and $G_n(\Lambda^+)$ is concave as a function of $\Lambda^+ \in (0, \infty)^{|\mathcal{I}_+(n)|}$ with a strictly negative definite (diagonal) Hessian matrix of second partial derivatives at each point, it follows from (35) that $K(n)$ is nonpositive and that it is strictly negative unless $\Lambda^+(n) = \Lambda^{+,*}(n)$. More precisely, by (35) and Taylor's theorem with remainder, for $v = \Lambda^+(n) - \Lambda^{+,*}(n)$,

$$\begin{aligned}-K(n) &= \nabla G_n(\Lambda^{+,*}(n)) \cdot v \\ &= G_n(\Lambda^+(n)) - G_n(\Lambda^{+,*}(n)) - \tfrac{1}{2} v \cdot (\nabla^2 G_n)(\tilde{\Lambda}) v,\end{aligned}$$

for some $\tilde{\Lambda}$ lying on the line segment between $\Lambda^+(n)$ and $\Lambda^{+,*}(n)$ and where $(\nabla^2 G_n)(\cdot)$ denotes the Hessian matrix for $G_n(\cdot)$. By the optimality of $\Lambda^+(n)$, $G_n(\Lambda^+(n)) \geq G_n(\Lambda^{+,*}(n))$, and by the strict negative definiteness of $(\nabla^2 G_n)(\tilde{\Lambda})$, it follows that the last line above is nonnegative and it is strictly positive unless $v = 0$. By Lemma 5.1, $v \equiv \Lambda^+(n) - \Lambda^{+,*}(n) = 0$ if and only if $n$ is an invariant state. □

COROLLARY 6.1. *At any regular point $t$ for a fluid model solution $n(\cdot)$, we have*

(37) $$\frac{d}{dt} F(n(t)) = K(n(t)) \leq 0,$$

*where the inequality is strict unless $n(t) \in \mathcal{M}_\alpha$.*

PROOF. This follows from Lemma 6.1 and (30)–(33). □

For each $n \in \mathbb{R}_+^{\mathbf{I}}$, let $w(n) = (w_j(n) : j \in \mathcal{J}_*)$ be defined by (20).

LEMMA 6.2. *For any fluid model solution $n(\cdot)$, $t \to w_j(n(t))$ is a nondecreasing function of $t \in [0, \infty)$ for each $j \in \mathcal{J}_*$.*

PROOF. Consider a fluid model solution $n(\cdot)$. Since $n(\cdot)$ is absolutely continuous, then so is the linear function $w(n(\cdot))$ of $n(\cdot)$. From (13) and (14) satisfied by a fluid model solution, at a regular point $t$ for $n(\cdot)$, we have for each $j \in \mathcal{J}_*$,

$$\text{(38)} \qquad \frac{d}{dt} w_j(n(t)) = \sum_{i \in \mathcal{I}_+(n(t))} A_{ji}(\rho_i - \Lambda_i(n(t))),$$

and (14) holds. Now, for $j \in \mathcal{J}_*$,

$$\text{(39)} \qquad C_j = \sum_{i \in \mathcal{I}} A_{ji}\rho_i,$$

and on substituting this into (14) we obtain

$$\text{(40)} \qquad \sum_{i \in \mathcal{I}_+(n(t))} A_{ji}\Lambda_i(n(t)) \leq \sum_{i \in \mathcal{I}_+(n(t))} A_{ji}\rho_i \qquad \text{for all } j \in \mathcal{J}_*,$$

and when combined with (38) this yields

$$\text{(41)} \qquad \frac{d}{dt} w_j(n(t)) \geq 0 \qquad \text{for all } j \in \mathcal{J}_*.$$

Since $w_j(n(\cdot))$ is obtained by integrating its almost everywhere defined derivative, it follows that $w_j(n(\cdot))$ is a nondecreasing function for each $j \in \mathcal{J}_*$. □

For each $w \in \mathbb{R}_+^{\mathbf{J}_*}$, define $\underline{F}(w)$ to be the optimal value attained in the optimization problem (22) and recall, from Section 5, the definition of $\Delta(w)$ as the optimizing value of $n$.

LEMMA 6.3. *The functions $\underline{F} : \mathbb{R}_+^{\mathbf{J}_*} \to \mathbb{R}_+$ and $\Delta : \mathbb{R}_+^{\mathbf{J}_*} \to \mathbb{R}_+^{\mathbf{I}}$ are continuous. In addition, $\underline{F}$ is a nondecreasing function, that is, if $w$ and $\tilde{w}$ are two vectors in $\mathbb{R}_+^{\mathbf{J}_*}$ such that $w \leq \tilde{w}$, then $\underline{F}(w) \leq \underline{F}(\tilde{w})$.*

PROOF. To prove the continuity of the functions $\underline{F}$ and $\Delta$, consider a sequence $\{w_m, m = 1, 2, \ldots\}$ in $\mathbb{R}_+^{\mathbf{J}_*}$ converging to some $w \in \mathbb{R}_+^{\mathbf{J}_*}$. By the growth property of $F$ that $F(n) \to \infty$ as $\|n\| \to \infty$, and the full row rank and nonnegativity assumptions on $A$, $\{\Delta(w_m), m = 1, 2, \ldots\}$ is bounded. Let $\tilde{n}$ be any cluster point of this sequence. Then, $\tilde{n}$ satisfies the constraints of the optimization problem (22) and, by the continuity of $F$, $F(\Delta(w_{m_k})) \to F(\tilde{n})$ as $k \to \infty$ for some subsequence $\{m_k\}$ of $\{m\}$. By the feasibility of $\tilde{n}$, $F(\tilde{n}) \geq F(\Delta(w))$. We claim that $\tilde{n} = \Delta(w)$.

For a proof by contradiction, suppose that $\tilde{n} \neq \Delta(w)$. Then, by the strict convexity of $F$, $\varepsilon \equiv F(\tilde{n}) - F(\Delta(w)) > 0$, and by the continuity of $F$, there is $\delta > 0$ such that $F(n) < F(\Delta(w)) + \varepsilon = F(\tilde{n})$ for all $n \in \mathbb{R}_+^{\mathbf{I}}$ satisfying $\|n - \Delta(w)\| < \delta$.

16    F. P. KELLY AND R. J. WILLIAMSSince $w_m \to w$ as $m \to \infty$ and $A$ has full row rank, there is $\hat{n} \in \mathbb{R}_+^{\mathbf{I}}$ such that $\|\hat{n} - \Delta(w)\| < \delta$ and $A_* M^{-1}\hat{n} \geq w_m$ for all $m$ sufficiently large. [Here $A_*$ is the $\mathbf{J}_* \times \mathbf{I}$ matrix of rank $\mathbf{J}_*$ obtained from $A$ by eliminating those rows of $A$ that are not indexed by $\mathcal{J}_*$ and $M = \operatorname{diag}(\mu)$.] To see this, let $A^\dagger$ be a $\mathbf{I} \times \mathbf{J}_*$ matrix that is a right inverse for $A_* M^{-1}$, that is, $A_* M^{-1} A^\dagger = I$, where $I$ denotes the $\mathbf{J}_* \times \mathbf{J}_*$ identity matrix. By the definition of $\Delta(w)$, $A_* M^{-1} \Delta(w) \geq w$ and so, since $w_m \to w$ as $m \to \infty$, there is $m(\delta)$ sufficiently large that for all $m \geq m(\delta)$,

$$A_* M^{-1} \Delta(w) \geq w_m - \frac{\delta}{2\sqrt{\mathbf{J}_*}\|A^\dagger\|}\mathbf{1},$$

where $\mathbf{1}$ denotes the $\mathbf{J}_*$-dimensional vector whose components are all 1's. Let $u = \frac{\delta}{2\sqrt{\mathbf{J}_*}\|A^\dagger\|} A^\dagger \mathbf{1}$. Note that $\|u\| \leq \frac{\delta}{2}$ and $A_* M^{-1} u = \frac{\delta}{2\sqrt{\mathbf{J}_*}\|A^\dagger\|}\mathbf{1}$. Define $\hat{n} = \Delta(w) + \tilde{u}$, where $\tilde{u}_i = |u_i|$ for $i \in \mathcal{I}$. Then, $\hat{n} \in \mathbb{R}_+^{\mathbf{I}}$, $\|\hat{n} - \Delta(w)\| = \|\tilde{u}\| = \|u\| < \delta$ and, since $A_* M^{-1}$ has nonnegative entries,

$$A_* M^{-1} \hat{n} \geq A_* M^{-1} \Delta(w) + A_* M^{-1} u \geq w_m$$

for all $m \geq m(\delta)$, as desired.

For $m \geq m(\delta)$, $\hat{n}$ is feasible for the optimization problem (22) with $w_m$ in place of $w$, and since $\Delta(w_m)$ is optimal for this problem, it follows that $F(\hat{n}) \geq F(\Delta(w_m))$. Hence, $F(\hat{n}) \geq \lim_k F(\Delta(w_{m_k})) = F(\tilde{n})$. But this yields a contradiction, since $F(\hat{n}) < F(\tilde{n})$ as $\|\hat{n} - \Delta(w)\| < \delta$. Thus, $\tilde{n} = \Delta(w)$, and since $\tilde{n}$ was an arbitrary cluster point of $\{\Delta(w_m), m = 1, 2, \dots\}$, it follows that this sequence converges to $\Delta(w)$, and

$$\underline{F}(w) = F(\Delta(w)) = \lim_m F(\Delta(w_m)) = \lim_m \underline{F}(w_m).$$

This implies the continuity of $\Delta$ and $\underline{F}$.

The nondecreasing property of $\underline{F}$ follows from the fact that, for $w$ and $\tilde{w}$ in $\mathbb{R}_+^{\mathbf{J}_*}$ satisfying $w \leq \tilde{w}$, any feasible solution for the problem (22) with $\tilde{w}$ in place of $w$ is feasible for the original problem with $w$.  □

We have the following characterization of the optimal solutions of (22).

LEMMA 6.4. *For each $w \in \mathbb{R}_+^{\mathbf{J}_*}$, a vector $n \in \mathbb{R}_+^{\mathbf{I}}$ is the unique optimal solution of (22) if and only if there is $p \in \mathbb{R}_+^{\mathbf{J}_*}$ such that for each $i \in \mathcal{I}$,*

$$(42) \qquad n_i = \rho_i \left( \frac{\sum_{j \in \mathcal{J}_*} p_j A_{ji}}{\kappa_i} \right)^{1/\alpha}$$

*and for each $j \in \mathcal{J}_*$,*

$$(43) \qquad p_j \left( \sum_{i \in \mathcal{I}} A_{ji} \frac{n_i}{\mu_i} - w_j \right) = 0$$

*and*

$$\sum_{i \in \mathcal{I}} A_{ji} \frac{n_i}{\mu_i} \geq w_j. \tag{44}$$

PROOF. Fix $w \in \mathbb{R}_+^{\mathbf{J}_*}$. For $n \in \mathbb{R}_+^{\mathbf{I}}$ and $p \in \mathbb{R}_+^{\mathbf{J}_*}$, let

$$L(n, p) = F(n) + \sum_{j \in \mathcal{J}_*} p_j \left( w_j - \sum_{i \in \mathcal{I}} A_{ji} \frac{n_i}{\mu_i} \right). \tag{45}$$

Suppose that $n \in \mathbb{R}_+^{\mathbf{I}}$ is the unique optimal solution of (22). Since $A$ has full row rank, no row of $A$ contains all zeros. Recall that $A$ only contains zeros and ones. It follows that there is $v \in \mathbb{R}_+^{\mathbf{I}}$ such that $\sum_{i \in \mathcal{I}} A_{ji} \frac{v_i}{\mu_i} > w_j$ for each $j \in \mathbf{J}_*$. Thus, Slater's constraint qualification (cf. [8], page 236) is satisfied and by the necessity theory of Lagrange multipliers (cf. [15], Theorem 1, page 217), there is $p \in \mathbb{R}_+^{\mathbf{J}_*}$ such that $n$ minimizes $L(\cdot, p)$ over $\mathbb{R}_+^{\mathbf{I}}$ and (43) holds for all $j \in \mathcal{J}_*$. Now, $L(\cdot, p)$ is continuously differentiable and so for each $i \in \mathcal{I}$ we have

$$\frac{\partial L}{\partial n_i}(n, p) \geq 0, \tag{46}$$

where the inequality is an equality when $n_i > 0$. If $n_i > 0$, this yields (42) immediately. If $n_i = 0$, this yields

$$0 = n_i \geq \rho_i \left( \frac{\sum_{j \in \mathcal{J}_*} p_j A_{ji}}{\kappa_i} \right)^{1/\alpha}. \tag{47}$$

Since $p$ has nonnegative components, it follows that the above holds with equality. Thus, (42) holds for all $i \in \mathcal{I}$. Inequality (44) holds for each $j \in \mathcal{J}_*$, since $n$ is feasible for (22).

Conversely, suppose that $n \in \mathbb{R}_+^{\mathbf{I}}$ and $p \in \mathbb{R}_+^{\mathbf{J}_*}$ satisfy (42)–(44). Then

$$\frac{\partial L}{\partial n_i}(n, p) = 0 \qquad \text{for all } i \in \mathcal{I}, \tag{48}$$

and it follows from the strict convexity of $F$ that $n$ is a global minimum for $L(\cdot, p)$ over $\mathbb{R}_+^{\mathbf{I}}$. By (44), $n$ is feasible for (22) and for any other feasible $\tilde{n} \in \mathbb{R}_+^{\mathbf{I}}$ we have

$$F(n) = L(n, p) \leq L(\tilde{n}, p) \leq F(\tilde{n}), \tag{49}$$

where we have used (43), the feasibility of $\tilde{n}$ and the fact that $p \in \mathbb{R}_+^{\mathbf{J}_*}$. Thus, $n$ is optimal for (22). $\square$

We can now prove Theorem 5.3.

PROOF OF THEOREM 5.3. By Theorem 5.1, $n \in \mathbb{R}_+^{\mathbf{I}}$ is an invariant state if and only if there is $p = q \in \mathbb{R}_+^{\mathbf{J}*}$ such that (42) holds for all $i \in \mathcal{I}$. Note that if $w = w(n)$, then (43) and (44) automatically hold for all $j \in \mathcal{J}_*$. Thus, on combining the above with Lemma 6.4, we see that $n \in \mathbb{R}_+^{\mathbf{I}}$ is an invariant state if and only if $n$ is the unique optimal of (22) with $w = w(n)$. The latter is equivalent to $n = \Delta(w(n))$. □

For the proof of the convergence result, Theorem 5.2, we introduce the following function $H$ and prove some of its elementary properties. For each $n \in \mathbb{R}_+^{\mathbf{I}}$, let

(50) $$H(n) = F(n) - \underline{F}(w(n)).$$

LEMMA 6.5. *The function $H : \mathbb{R}_+^{\mathbf{I}} \to \mathbb{R}$ is continuous. Furthermore, it is zero on the set of invariant states, $\mathcal{M}_\alpha$, and it is strictly positive on $\mathbb{R}_+^{\mathbf{I}} \setminus \mathcal{M}_\alpha$.*

PROOF. The continuity of $H$ on $\mathbb{R}_+^{\mathbf{I}}$ follows from the facts that $F$ is clearly continuous, $\underline{F}$ is continuous by Lemma 6.3 and the linear function $w(\cdot)$ is continuous.

Since $n \in \mathbb{R}_+^{\mathbf{I}}$ is feasible for (22) with $w = w(n)$, we have $F(n) \geq \underline{F}(w(n))$ and hence $H(n) \geq 0$. Furthermore, $H(n) = 0$ if and only if $n$ is the optimal solution for (22) with $w = w(n)$, that is, if and only if $n = \Delta(w(n))$. The latter occurs if and only if $n$ is an invariant state by Theorem 5.3. □

LEMMA 6.6. *For any fluid model solution $n(\cdot)$, $t \to H(n(t))$ is a nonincreasing function of $t \in [0, \infty)$.*

PROOF. By Corollary 6.1, $F(n(t)) = F(n(0)) + \int_0^t K(n(s))\,ds$ is a nonincreasing function of $t$, and, by the combination of Lemmas 6.2 and 6.3, $\underline{F}(w(n(t)))$ is a nondecreasing function of $t$, and so it follows that $H(n(t)) = F(n(t)) - \underline{F}(w(n(t)))$ is a nonincreasing function of $t$. □

REMARK 6.1. A stronger form of Lemma 6.6, which shows that $H(n(t))$ is strictly decreasing at times where $n(t) \notin \mathcal{M}_\alpha$ will be developed and used in the proof of Theorem 5.2 given below.

PROOF OF THEOREM 5.2. Fix $R \in (0, \infty)$ and $\varepsilon > 0$. Let

$$\hat{F}(R) = \sup\{F(v) : v \in \mathbb{R}_+^{\mathbf{I}}, \|v\| \leq R\}.$$

Since $F$ is continuous on $\mathbb{R}_+^{\mathbf{I}}$, $\hat{F}(R)$ is finite. For any fluid model solution $n(\cdot)$ satisfying $\|n(0)\| \leq R$, on integrating (37) we see that $F(n(t)) \leq F(n(0))$ for all $t \geq 0$. The fact that $F(n) \to \infty$ as $\|n\| \to \infty$ implies that there is a closed ball $B$ (depending on $R$) in $\mathbb{R}_+^{\mathbf{I}}$ that is centered at the origin and of finite radius such that $F(v) > \hat{F}(R)$ for all $v \in B^c$, where $B^c$ denotes the complement of $B$ in $\mathbb{R}_+^{\mathbf{I}}$. Combining the above, we see that for any fluid model solution $n(\cdot)$ satisfying $\|n(0)\| \leq R$, we have $n(t) \in B$ for all $t \in [0, \infty)$.

Let

$$\tag{51} D = \{v \in B : d(\mathcal{M}_\alpha, v) \geq \varepsilon\}.$$

Note that $D$ depends on $R$ and $\varepsilon$. By Lemma 6.5, the function $H$ is continuous and strictly positive on the compact set $D$. It follows that $\delta \equiv \inf\{H(v) : v \in D\}$ is strictly positive, and the set $\tilde{D} \equiv \{v \in B : H(v) \geq \delta\}$ contains $D$. Moreover, since $H$ is zero on $\mathcal{M}_\alpha$, $\tilde{D}$ does not meet $\mathcal{M}_\alpha$.

Consider a fluid model solution $n(\cdot)$ satisfying $\|n(0)\| \leq R$. Let $\tilde{T}(n) = \inf\{t \geq 0 : n(t) \notin \tilde{D}\}$. Since $n(\cdot)$ remains in $B$ for all time, it follows that if $\tilde{T}(n) < \infty$, then $n(\cdot)$ exits $\tilde{D}$ by violating the constraint $H(n(\cdot)) \geq \delta$. Then, since $H(n(\cdot))$ is a nonincreasing function (by Lemma 6.6), it follows that $n(t) \notin \tilde{D}$ for all $t > \tilde{T}(n)$. Consequently, since $D \subset \tilde{D}$, we have $d(\mathcal{M}_\alpha, n(t)) < \varepsilon$ for all $t > \tilde{T}(n)$. We now develop an upper bound on $\tilde{T}(n)$. For $0 \leq t \leq \tilde{T}(n)$,

$$\tag{52} H(n(t)) - H(n(0)) \leq F(n(t)) - F(n(0))$$

$$\tag{53} = \int_0^t K(n(s))\, ds,$$

where we have used the nondecreasing property of $\underline{F}(w(n(\cdot)))$ for the first line, and Corollary 6.1 for the second line. By Lemma 6.1, $K$ is continuous and strictly negative off the manifold $\mathcal{M}_\alpha$. It follows that there is $C_{R,\varepsilon} > 0$ such that $K$ is bounded above by $-C_{R,\varepsilon}$ on the compact set $\tilde{D}$ which does not meet $\mathcal{M}_\alpha$. Then, the above yields

$$\tag{54} \tilde{T}(n) \leq \bar{H}_R / C_{R,\varepsilon},$$

where $\bar{H}_R = \sup\{H(v) : \|v\| \leq R\} < \infty$, and the desired result follows. $\square$

## APPENDIX A:

**Properties of the bandwidth sharing policy.** For each $n \in \mathbb{R}_+^{\mathbf{I}} \setminus \{0\}$ and $\Lambda^+ \in \mathbb{R}_+^{|\mathcal{I}_+(n)|}$, recall the definition of $G_n(\Lambda^+)$ from (4).

LEMMA A.1. *For each $n \in \mathbb{R}_+^{\mathbf{I}} \setminus \{0\}$, there is a unique optimal solution $\Lambda^+(n) = (\Lambda_i : i \in \mathcal{I}_+(n))$ of (1)–(3); each of its components is strictly positive and uniformly bounded by $C^* = \max_{j \in \mathcal{J}} C_j$.*

PROOF. Fix $n \in \mathbb{R}_+^{\mathbf{I}} \setminus \{0\}$. If $\alpha \in (0,1)$, the objective function $G_n(\Lambda^+)$ is continuous and strictly concave as a function of $\Lambda^+ = (\Lambda_i : i \in \mathcal{I}_+(n))$ on

the compact set

$$(55) \quad \left\{\Lambda^+ : \Lambda_i \geq 0 \text{ for all } i \in \mathcal{I}_+(n), \sum_{i \in \mathcal{I}_+(n)} A_{ji}\Lambda_i \leq C_j \text{ for all } j \in \mathcal{J}\right\}.$$

It follows that the maximum of the objective function is achieved on this set, and by the strict concavity, the maximizing point is unique. Furthermore, since the derivative of the objective function with respect to $\Lambda_i$ diverges to $+\infty$ as $\Lambda_i \downarrow 0$ for each $i \in \mathcal{I}_+(n)$, the maximizing point must have strictly positive components, that is, the maximum cannot occur with one of the $\Lambda_i$, $i \in \mathcal{I}_+(n)$, equal to zero.

For $\alpha \in [1,\infty)$, by convention, the objective function takes the value $-\infty$ whenever $\Lambda_i = 0$ for some $i \in \mathcal{I}_+(n)$. Thus, any optimal solution to (1)–(3) has all components strictly positive. In fact, if we let $c = \min_{j \in \mathcal{J}} C_j/|\mathcal{I}|$, then the $|\mathcal{I}_+(n)|$-dimensional vector $\tilde{\Lambda}^+$ with all components set equal to $c$ satisfies the constraints in (1)–(3). By combining the fact that the feasible set for (1)–(3) is bounded with the fact that the objective function diverges to $-\infty$ as $\Lambda_i \to 0$ for any $i \in \mathcal{I}_+(n)$, we may conclude that there is $\delta > 0$ such that for all $\Lambda^+ = (\Lambda_i : i \in \mathcal{I}_+(n))$ satisfying the constraints in (1)–(3) and such that $\Lambda_i < \delta$ for at least one $i \in \mathcal{I}_+(n)$, we have $G_n(\Lambda^+) < G_n(\tilde{\Lambda}^+)$. It follows that any optimal solution of (1)–(3) must be in the compact set

$$(56) \quad \left\{\Lambda^+ : \delta \leq \Lambda_i \text{ for all } i \in \mathcal{I}_+(n), \sum_{i \in \mathcal{I}_+(n)} A_{ji}\Lambda_i \leq C_j \text{ for all } j \in \mathcal{J}\right\}.$$

Since the objective function is continuous and strictly concave on this compact set, it has a unique maximum there, which is the optimal solution of (1)–(3).

Since for each $i \in \mathcal{I}$, there is $j \in \mathcal{J}$, such that $A_{ji} = 1$, it follows that $\Lambda_i^+(n) \leq C^* = \max_{j \in \mathcal{J}} C_j$ for all $i \in \mathcal{I}_+(n)$.  $\square$

In the following lemmas, $\Lambda : \mathbb{R}_+^{\mathbf{I}} \to \mathbb{R}_+^{\mathbf{I}}$ is as defined in Section 3.

LEMMA A.2. *Fix $n \in \mathbb{R}_+^{\mathbf{I}}$. Then, $\Lambda(rn) = \Lambda(n)$ for all $r > 0$.*

PROOF. Fix $r > 0$. Note that $\mathcal{I}_+(n) = \mathcal{I}_+(rn)$ and so $\Lambda_i(n) = \Lambda_i(rn) = 0$ for all $i \in \mathcal{I}_0(n) = \mathcal{I} \setminus \mathcal{I}_+(n)$. Also, for $\Lambda^+ \in \mathbb{R}_+^{|\mathcal{I}_+(n)|}$, $G_{rn}(\Lambda^+) = r^\alpha G_n(\Lambda^+)$ [for all $\alpha \in (0,\infty)$], and since the constraints in (1)–(3) do not involve $n$, it follows that the optimal solution $\Lambda^+(rn)$ for the objective function $G_{rn}(\cdot)$ is the same as the optimal solution $\Lambda^+(n)$ for the objective function $G_n(\cdot)$, and hence $\Lambda^+(rn) = \Lambda^+(n)$.  $\square$

LEMMA A.3. *For each $i \in \mathcal{I}$, the $i$th component $\Lambda_i(\cdot)$ of the function $\Lambda(\cdot)$ is continuous on $\{n \in \mathbb{R}_+^{\mathbf{I}} : n_i > 0\}$.*

PROOF. We first prove continuity of $\Lambda(\cdot)$ at $n^* \in \mathbb{R}_+^{\mathbf{I}}$ satisfying $n_i^* > 0$ for all $i \in \mathcal{I}$. In this case, $\mathcal{I}_+(n^*) = \mathcal{I}$ and, by Lemma A.1, the optimal solution $\Lambda(n^*)$ to the problem (1)–(3) with $n = n^*$ has all components strictly positive. Given $\varepsilon > 0$, let $B_\varepsilon$ be a nonempty open ball that is centered at $\Lambda(n^*)$, that has radius less than or equal to $\varepsilon$ and that is a positive distance from the boundary of the orthant $\mathbb{R}_+^{\mathbf{I}}$. Let $D_\varepsilon$ denote the compact set of $\Lambda \in \mathbb{R}_+^{\mathbf{I}}$ that satisfy the constraints of (1)–(3) with $n = n^*$ there and that are outside of $B_\varepsilon$.

We claim that there is $\eta > 0$, $\beta \in (0, \min_{i=1}^{\mathbf{I}} n_i^*)$ and $\delta > 0$, such that

(57) $$G_{n^*}(\Lambda) < G_{n^*}(\Lambda(n^*)) - \eta \qquad \text{for all } \Lambda \in D_\varepsilon$$

and

(58) $$G_n(\Lambda) < G_{n^*}(\Lambda(n^*)) - \eta/2,$$

whenever $n \in \mathbb{R}_+^{\mathbf{I}}$ satisfies $\|n - n^*\| < \beta$ and $\Lambda \in D_\varepsilon$ is such that $\Lambda_i < \delta$ for at least one $i \in \mathcal{I}$. (Note that $n_i > 0$ for all $i$ when $\|n - n^*\| < \beta$.) The first inequality (57) can be proved by contradiction. For if (57) does not hold for some $\eta > 0$, then for each positive integer $k$, there is $\Lambda^k \in D_\varepsilon$ such that

(59) $$G_{n^*}(\Lambda^k) \geq G_{n^*}(\Lambda(n^*)) - 1/k.$$

By the compactness of $D_\varepsilon$, there is $\tilde{\Lambda} \in D_\varepsilon$ such that $\{\Lambda^k\}$ converges to $\tilde{\Lambda}$ along some subsequence. If $\alpha \in (0, 1)$, then by the continuity of $G_{n^*}(\cdot)$ on $[0, \infty)^{\mathbf{I}}$ in this case, it follows on passing to the limit in (59) that $G_{n^*}(\tilde{\Lambda}) \geq G_{n^*}(\Lambda(n^*))$, which contradicts the uniqueness of the optimum $\Lambda(n^*)$. If $\alpha \in [1, \infty)$ and $\tilde{\Lambda} \in (0, \infty)^{\mathbf{I}}$, then the continuity of $G_{n^*}(\cdot)$ on $(0, \infty)^{\mathbf{I}}$ yields a contradiction in the same manner as for $\alpha \in (0, 1)$. Finally, if $\alpha \in [1, \infty)$ and $\tilde{\Lambda}_i = 0$ for some $i \in \mathcal{I}$, then since $\Lambda_i^k \to \tilde{\Lambda}_i$ along some subsequence of $k$'s, either $\kappa_i n_i^* \log \Lambda_i^k \to -\infty$ (if $\alpha = 1$) or $\kappa_i (n_i^*)^\alpha (\Lambda_i^k)^{1-\alpha}/(1-\alpha) \to -\infty$ (if $\alpha > 1$), along this same subsequence. The other terms (indexed by $l \neq i$) in the sum constituting $G_{n^*}(\Lambda^k)$ are either bounded (if $\tilde{\Lambda}_l \neq 0$) or go to $-\infty$ (if $\tilde{\Lambda}_l = 0$) as $k \to \infty$ along the subsequence. Consequently, $G_{n^*}(\Lambda^k) \to -\infty$ as $k \to \infty$ along the subsequence, which contradicts (59). Now we show the result pertaining to (58). For $\alpha \in (0, 1)$, it follows from the uniform continuity of $G_n(\Lambda)$ as a function of $(n, \Lambda)$ on any compact subset of $[0, \infty)^{\mathbf{I}} \times [0, \infty)^{\mathbf{I}}$ that there is $\beta \in (0, \min_{i=1}^{\mathbf{I}} n_i^*)$ such that

$$G_n(\Lambda) < G_{n^*}(\Lambda) + \eta/2$$

for all $n \in \mathbb{R}_+^{\mathbf{I}}$ satisfying $\|n - n^*\| \leq \beta$ and $\Lambda \in D_\varepsilon$. Inequality (58) follows immediately from this together with (57). For $\alpha = 1$ and fixed $\beta \in (0, \min_{l=1}^{\mathbf{I}} n_l^*)$, if $n \in \mathbb{R}_+^{\mathbf{I}}$ satisfies $\|n - n^*\| < \beta$, $\delta \in (0, 1)$, $\Lambda \in D_\varepsilon$ and $i \in \mathcal{I}$ such that $\Lambda_i < \delta$, then we have

$$G_n(\Lambda) = \sum_{l \in \mathcal{I}} \kappa_l n_l \log \Lambda_l \tag{60}$$

$$\leq \kappa_i(n_i^* - \beta) \log \delta + \sum_{l \neq i} \kappa_l(n_l^* + \beta) \log(C^* \vee 1), \tag{61}$$

where $C^* = \max_{j \in \mathcal{J}} C_j$. Since $n_i^* - \beta > 0$, by the definition of $\beta$, the last line above goes to $-\infty$ as $\delta \to 0$. As there are only finitely many indices $i \in \mathcal{I}$ and the right member of (58) does not vary with $n$ or $\Lambda$, it follows that for any fixed $\beta \in (0, \min_{l=1}^{\mathbf{I}} n_l^*)$, there is $\delta \in (0,1)$ such that whenever $n \in \mathbb{R}_+^{\mathbf{I}}$ satisfies $\|n - n^*\| < \beta$, $\Lambda \in D_\varepsilon$ and $\Lambda_i < \delta$ for some $i \in \mathcal{I}$, then (58) holds. A similar proof of (58) holds for $\alpha \in (1, \infty)$, with the exception that (61) is replaced with

$$G_n(\Lambda) \leq \kappa_i(n_i^* - \beta)^\alpha \frac{\delta^{1-\alpha}}{1-\alpha}. \tag{62}$$

Since $G_n(\Lambda)$ is uniformly continuous as a function of $(n, \Lambda)$ on any compact subset of $(0, \infty)^{\mathbf{I}} \times (0, \infty)^{\mathbf{I}}$, there is $\gamma \in (0, \beta/2)$ such that for all $n \in \mathbb{R}_+^{\mathbf{I}}$ satisfying $\|n - n^*\| \leq \gamma$ and $\Lambda \in D_\varepsilon$ satisfying $\Lambda_i \geq \delta$ for all $i \in \mathcal{I}$, we have

$$G_n(\Lambda) < G_{n^*}(\Lambda) + \eta/2 \tag{63}$$

and

$$G_n(\Lambda(n^*)) > G_{n^*}(\Lambda(n^*)) - \eta/2. \tag{64}$$

Combining (57) and (58) with (63) and (64), we see that for all $n \in \mathbb{R}_+^{\mathbf{I}}$ satisfying $\|n - n^*\| \leq \gamma$ and all $\Lambda \in D_\varepsilon$, we have

$$G_n(\Lambda) < G_{n^*}(\Lambda(n^*)) - \eta/2 < G_n(\Lambda(n^*)). \tag{65}$$

It follows that for all $n \in \mathbb{R}_+^{\mathbf{I}}$ satisfying $\|n - n^*\| \leq \gamma$, the optimal solution $\Lambda(n)$ of (1)–(3) must lie within $B_\varepsilon$, and hence must be no more than distance $\varepsilon$ from $\Lambda(n^*)$. Hence, $\Lambda(\cdot)$ is continuous at $n^*$.

Next, consider $n^* \in \mathbb{R}_+^{\mathbf{I}} \setminus \{0\}$ such that $n_k^* = 0$ for at least one $k$. We will show that as a function of $n$, $(\Lambda_i(n), i \in \mathcal{I}_+(n^*))$ is continuous at $n = n^*$. Small perturbations $n$ of $n^*$ that only involve changes to the strictly positive components of $n^*$ can be handled in a similar manner to that in the first three paragraphs of this proof. It is perturbations for which $\mathcal{I}_+(n) \neq \mathcal{I}_+(n^*)$ that require additional argument, which we give below. For fixed $\varepsilon > 0$, let $B_\varepsilon$ be a nonempty open ball in $\mathbb{R}_+^{|\mathcal{I}_+(n^*)|}$ that is centered at $\Lambda^+(n^*) = (\Lambda_i(n^*) : i \in \mathcal{I}_+(n^*))$, that has radius less than or equal to $\varepsilon$ and that is a positive distance from the boundary of the orthant $\mathbb{R}_+^{|\mathcal{I}_+(n^*)|}$. Let $D_\varepsilon$ denote the set of $\tilde{\Lambda} \in \mathbb{R}_+^{|\mathcal{I}_+(n^*)|}$ that satisfy the constraints of (1)–(3) with $n = n^*$ there and that are outside of $B_\varepsilon$. In a similar manner to that for (57), by the strict concavity of $G_{n^*}(\cdot)$ on the interior of $\mathbb{R}_+^{|\mathcal{I}_+(n^*)|}$ and its behavior near the boundary of $\mathbb{R}_+^{|\mathcal{I}_+(n^*)|}$, there is $\eta > 0$ such that

$$(66) \qquad G_{n^*}(\tilde{\Lambda}) < G_{n^*}(\Lambda^+(n^*)) - \eta \qquad \text{for all } \tilde{\Lambda} \in D_\varepsilon.$$

We now need separate arguments for the cases $\alpha < 1, \alpha > 1$ and $\alpha = 1$.

We first consider $\alpha \in (0,1)$. For a vector $n \in \mathbb{R}_+^{\mathbf{I}}$ such that $n_i > 0$ for all $i \in \mathcal{I}_+(n^*)$, the vector $\Lambda^+(n^*) \in \mathbb{R}_+^{|\mathcal{I}_+(n^*)|}$ can be extended, by the addition of zero components, to a feasible vector in $\mathbb{R}_+^{|\mathcal{I}_+(n)|}$ for the optimization problem (1)–(3) associated with $n$. Using the fact that

$$(67) \qquad \kappa_i n_i^\alpha \frac{\Lambda_i^{1-\alpha}}{1-\alpha}$$

is zero when $\Lambda_i = 0$ and $\alpha \in (0,1)$, it follows that the optimal solution $\Lambda^+(n)$ of the optimization problem (1)–(3) satisfies

$$(68) \qquad G_n(\Lambda^+(n)) \geq G_{\tilde{n}}(\Lambda^+(n^*)),$$

where $\tilde{n}$ is the vector in $\mathbb{R}_+^{\mathbf{I}}$ obtained by resetting the entries in $n$ indexed by $i \in \mathcal{I}_+(n) \setminus \mathcal{I}_+(n^*)$ to zero. Let $\tilde{\Lambda}^+(n)$ denote the vector in $\mathbb{R}_+^{|\mathcal{I}_+(n^*)|}$ consisting of the components of $\Lambda^+(n)$ indexed by $i \in \mathcal{I}_+(n^*)$. Then,

$$(69) \qquad G_n(\Lambda^+(n)) = G_{\tilde{n}}(\tilde{\Lambda}^+(n)) + \sum_{i \in \mathcal{I}_+(n) \setminus \mathcal{I}_+(n^*)} \kappa_i n_i^\alpha \frac{(\Lambda_i^+(n))^{1-\alpha}}{1-\alpha}.$$

Since (67) tends to zero as $n_i \to 0$, uniformly for $\Lambda_i \in [0, C^*]$, there is $\gamma_1 > 0$ such that whenever $\|n - n^*\| < \gamma_1$, we have $n_i \in (0, \infty)$ for all $i \in \mathcal{I}_+(n^*)$ and the sum on the right-hand side of (69) is less than $\eta/6$. Combining this with (68), we see that for $\|n - n^*\| < \gamma_1$,

$$(70) \qquad G_{\tilde{n}}(\tilde{\Lambda}^+(n)) \geq G_{\tilde{n}}(\Lambda^+(n^*)) - \eta/6.$$

By the uniform continuity of $G_{\tilde{n}}(\Lambda^+(n^*)) = \sum_{i \in \mathcal{I}_+(n^*)} \kappa_i n_i^\alpha \frac{(\Lambda_i^+(n^*))^{1-\alpha}}{1-\alpha}$ as a function of $(n_i : i \in \mathcal{I}_+(n^*))$ on compact subsets of $(0, \infty)^{|\mathcal{I}_+(n^*)|}$, there is $\gamma_2 \in (0, \gamma_1)$ such that

$$(71) \qquad \|G_{\tilde{n}}(\Lambda^+(n^*)) - G_{n^*}(\Lambda^+(n^*))\| < \eta/6,$$

whenever $\|n - n^*\| < \gamma_2$. Combining (70) with (71), we obtain

$$(72) \qquad G_{\tilde{n}}(\tilde{\Lambda}^+(n)) \geq G_{n^*}(\Lambda^+(n^*)) - \eta/3,$$

for all $n \in \mathbb{R}_+^{\mathbf{I}}$ satisfying $\|n - n^*\| < \gamma_2$. Finally, by the same kind of argument as in the second and third paragraphs of this proof [cf. (58)–(65)], but with $\Lambda$ replaced by $(\Lambda_i : i \in \mathcal{I}_+(n^*))$ and $n$ replaced by $(n_i : i \in \mathcal{I}_+(n^*))$, using (66), there is $\gamma_3 \in (0, \gamma_2)$ such that for all $n \in \mathbb{R}_+^{\mathbf{I}}$ satisfying $\|n - n^*\| < \gamma_3$ and $\tilde{\Lambda} \in D_\varepsilon$, we have

$$G_{n^*}(\Lambda^+(n^*)) > G_{\tilde{n}}(\tilde{\Lambda}) + \eta/2. \tag{73}$$

Combining the above, we have that for all $n \in \mathbb{R}_+^{\mathbf{I}}$ such that $\|n - n^*\| < \gamma_3$ and $\tilde{\Lambda} \in D_\varepsilon$,

$$G_{\tilde{n}}(\tilde{\Lambda}^+(n)) > G_{\tilde{n}}(\tilde{\Lambda}) + \eta/6. \tag{74}$$

It then follows that for $\|n - n^*\| < \gamma_3$, $\tilde{\Lambda}^+(n)$ is not in $D_\varepsilon$ and it satisfies the constraints of the optimization problem (1)–(3) with $\tilde{n}$ (or $n^*$) in place of $n$, and so it must be in $B_\varepsilon$. This completes the proof of the desired continuity for the case $\alpha \in (0, 1)$.

For the case $\alpha \in (1, \infty)$, for $n_i > 0$, the term (67) is negative, and is unbounded below as $\Lambda_i \downarrow 0$. However the choice $\Lambda_i = n_i$ gives an evaluation $\kappa_i n_i/(1-\alpha)$ that converges to zero as $n_i \downarrow 0$. Since $\sum_{i \in \mathcal{I}_+(n^*)} \kappa_i n_i^\alpha \frac{\Lambda_i^{1-\alpha}}{1-\alpha}$ is uniformly continuous as a function of $((n_i, \Lambda_i): i \in \mathcal{I}_+(n^*))$ on compact subsets of $((0, \infty) \times (0, \infty))^{|\mathcal{I}_+(n^*)|}$, there is $\gamma_1 > 0$ such that whenever $\|n_i - n_i^*\| < \gamma_1$ and $\|\Lambda_i - \Lambda_i^+(n^*)\| < \gamma_1$ for all $i \in \mathcal{I}_+(n^*)$, we have $n_i > 0, \Lambda_i > 0$ for all $i \in \mathcal{I}_+(n^*)$ and

$$\left\| \sum_{i \in \mathcal{I}_+(n^*)} \kappa_i n_i^\alpha \frac{\Lambda_i^{1-\alpha}}{1-\alpha} - G_{n^*}(\Lambda^+(n^*)) \right\| < \frac{\eta}{6}. \tag{75}$$

Since $n_i^* = 0$ for all $i \notin \mathcal{I}_+(n^*)$, there is $\gamma_2 \in (0, \gamma_1)$ such that

$$-\frac{\eta}{6} < \frac{1}{1-\alpha} \sum_{i \in \mathcal{I}_+(n) \setminus \mathcal{I}_+(n^*)} \kappa_i n_i < 0 \tag{76}$$

and

$$\sum_{i \in \mathcal{I}_+(n) \setminus \mathcal{I}_+(n^*)} A_{ji} n_i < (\gamma_1/2) \wedge C_j \qquad \text{for all } j \in \mathcal{J}, \tag{77}$$

whenever $n \in \mathbb{R}_+^{\mathbf{I}}$ satisfies $\|n - n^*\| < \gamma_2$. Then, for such $n$, the vector $\breve{\Lambda} \in \mathbb{R}_+^{|\mathcal{I}_+(n)|}$ defined by $\breve{\Lambda}_i = n_i$ for $i \in \mathcal{I}_+(n) \setminus \mathcal{I}_+(n^*)$ and $\breve{\Lambda}_i = \Lambda_i^+(n^*) - \gamma_1/2$ for $i \in \mathcal{I}_+(n^*)$, satisfies the constraints of the optimization problem (1)–(3) for $n$. Since $\Lambda^+(n)$ is the optimal solution of that problem, using the same notation with tildes as for the case $\alpha \in (0, 1)$, it follows, by the negativity of (67), that

$$G_{\tilde{n}}(\tilde{\Lambda}^+(n)) \geq G_n(\Lambda^+(n)) \geq G_n(\breve{\Lambda}), \tag{78}$$

where, by (75) and (76),

$$(79) \quad G_n(\check{\Lambda}) = \frac{1}{1-\alpha} \sum_{i \in \mathcal{I}_+(n^*)} \kappa_i n_i^\alpha \check{\Lambda}_i^{1-\alpha} + \frac{1}{1-\alpha} \sum_{i \in \mathcal{I}_+(n) \setminus \mathcal{I}_+(n^*)} \kappa_i n_i$$

$$(80) \qquad \geq G_{n^*}(\Lambda^+(n^*)) - \frac{\eta}{6} - \frac{\eta}{6}.$$

The proof proceeds in the same manner as from (72) onward in the proof for $\alpha \in (0,1)$, to show that $\tilde{\Lambda}^+(n)$ is in $B_\varepsilon$ whenever $n$ is sufficiently close to $n^*$.

If $\alpha = 1$ and $n_i > 0$, the term (67) becomes $\kappa_i n_i \log \Lambda_i$, which may take positive or negative values. Nevertheless, one can proceed in a very similar manner to that for the case $\alpha \in (1, \infty)$, after observing that if $\Lambda_i = n_i$, then the evaluation $\kappa_i n_i \log n_i$ is negative for small enough $n_i$ and approaches zero as $n_i \downarrow 0$, while the maximum over $0 \leq \Lambda_i \leq C^*$ of the analogue of (67), namely $\kappa_i n_i \log C^*$, also approaches zero as $n_i \downarrow 0$. It follows that there is $\gamma_2 > 0$ such that for $n \in \mathbb{R}_+^\mathbf{I}$ satisfying $\|n - n^*\| < \gamma_2$, we have $n_i > 0$ for all $i \in \mathcal{I}_+(n^*)$ and

$$(81) \qquad G_{\tilde{n}}(\tilde{\Lambda}^+(n)) \geq G_n(\Lambda^+(n)) - \eta/12$$

$$(82) \qquad \qquad \geq G_{n^*}(\Lambda^+(n^*)) - \eta/12 - \eta/12 - \eta/6.$$

The proof can then be completed as in the case $\alpha \in (1, \infty)$. $\square$

REMARK A.1. $\Lambda_i(n)$ may be discontinuous at $n$ if $n_i = 0$.

For $n \in \mathbb{R}_+^\mathbf{I} \setminus \{0\}$, consider the Lagrangian

$$L_n(\Lambda^+, p) = G_n(\Lambda^+) + \sum_{j \in \mathcal{J}} p_j \left( C_j - \sum_{i \in \mathcal{I}_+(n)} A_{ji} \Lambda_i \right),$$

where $G_n(\Lambda^+)$ is defined by (4), $\Lambda^+ = (\Lambda_i : i \in \mathcal{I}_+(n)) \in (0, \infty)^{|\mathcal{I}_+(n)|}$ and $p = (p_j : j \in \mathcal{J}) \in \mathbb{R}_+^\mathbf{J}$ is a vector of Lagrange multipliers. Since the optimal solution $\Lambda^+(n)$ of (1)–(3) has strictly positive components (cf. Lemma A.1), we can use the theory of Lagrange multipliers to characterize $\Lambda^+(n)$ as follows.

LEMMA A.4. *Fix* $n \in \mathbb{R}_+^\mathbf{I} \setminus \{0\}$. *A vector* $\Lambda^+ = (\Lambda_i : i \in \mathcal{I}_+(n))$ *is the unique optimal solution of* (1)–(3) *if and only if there is* $p \in \mathbb{R}_+^\mathbf{J}$ *such that*

$$(83) \qquad p_j \left( C_j - \sum_{i \in \mathcal{I}_+(n)} A_{ji} \Lambda_i \right) = 0 \qquad \text{for all } j \in \mathcal{J},$$

$$(84) \qquad \sum_{j \in \mathcal{J}} p_j A_{ji} > 0 \qquad \text{for all } i \in \mathcal{I}_+(n),$$

$$\Lambda_i = n_i \left( \frac{\kappa_i}{\sum_{j \in \mathcal{J}} p_j A_{ji}} \right)^{1/\alpha} \qquad \text{for all } i \in \mathcal{I}_+(n) \tag{85}$$

*and*

$$\sum_{i \in \mathcal{I}_+(n)} A_{ji} \Lambda_i \leq C_j \qquad \text{for all } j \in \mathcal{J}. \tag{86}$$

PROOF. By Lemma A.1, all of the components of the unique optimal solution $\Lambda^+(n)$ of (1)–(3) are strictly positive. Since $C_j > 0$ for each $j \in \mathcal{J}$, there is $\Lambda^+ \in (0, \infty)^{|\mathcal{I}_+(n)|}$ satisfying $C_j - \sum_{i \in \mathcal{I}_+(n)} A_{ji} \Lambda_i > 0$ for each $j \in \mathcal{J}$. It follows that Slater's constraint qualification (cf. [8], page 236) is satisfied and so by the necessity theory of Lagrange multipliers (cf. [15], Theorem 1, page 217), there is $p \in \mathbb{R}_+^{\mathbf{J}}$ such that $\Lambda^+(n)$ maximizes $L_n(\Lambda^+, p)$ over $\Lambda^+ \in (0, \infty)^{|\mathcal{I}_+(n)|}$, and (83) holds with $\Lambda_i = \Lambda_i^+(n)$ for all $i \in \mathcal{I}_+(n)$. Since $\Lambda^+ = \Lambda^+(n)$ is an interior maximum for $L_n(\cdot, p)$, which is continuously differentiable as a function of $\Lambda^+ \in (0, \infty)^{|\mathcal{I}_+(n)|}$, it follows that

$$0 = \frac{\partial L_n}{\partial \Lambda_i}(\Lambda^+(n), p) = \kappa_i \left(\frac{n_i}{\Lambda_i(n)}\right)^\alpha - \sum_{j \in \mathcal{J}} p_j A_{ji} \qquad \text{for all } i \in \mathcal{I}_+(n). \tag{87}$$

The relations (84) and (85) follow immediately from this and the strict positivity of $\Lambda^+(n)$. Condition (86) follows from the constraints satisfied by $\Lambda^+(n)$.

Conversely, if $p \in \mathbb{R}_+^{\mathbf{J}}$ and $\Lambda^+ = (\Lambda_i : i \in \mathcal{I}_+(n))$ satisfy (83)–(86), then using the strict concavity of $L_n(\cdot, p)$ one can verify that $(\Lambda_+, p)$ is a saddle point of $L_n(\cdot, \cdot)$ on $(0, \infty)^{|\mathcal{I}_+(n)|} \times \mathbb{R}_+^{\mathbf{J}}$ and so by a sufficiency theorem for Lagrange multipliers (cf. [15], Theorem 2, page 221), $\Lambda^+$ is an optimal solution of (1)–(3) with $\Lambda_i > 0$ in place of $\Lambda_i \geq 0$ for all $i \in \mathcal{I}_+(n)$. But since the optimal solution of (1)–(3) has strictly positive components, it follows that $\Lambda^+$ is optimal for (1)–(3). □

REMARK A.2. Given $n \in \mathbb{R}_+^{\mathbf{I}} \setminus \{0\}$, while the vector $p = p(n)$ may be non-unique, the vector $(\sum_{j \in \mathcal{J}} p_j(n) A_{ji}, i \in \mathcal{I}_+(n))$ is unique by (85) and the fact that $\Lambda^+(n)$ is unique, by Lemma A.1.

## APPENDIX B:

**Fluid model solutions as fluid limits.** We shall use the following representation for the Markov process $N$ introduced in Section 4:

$$N_i(t) = N_i(0) + E_i(t) - S_i(T_i(t)), \qquad i \in \mathcal{I}, \ t \geq 0, \tag{88}$$

where $N(0)$ is a random variable taking values in $\mathbb{R}_+^{\mathbf{I}}$ that is independent of $E, S$; $E, S$ are independent $\mathbf{I}$-dimensional processes with independent components such that for each $i \in \mathcal{I}$, $E_i$ is a Poisson process with rate $\nu_i$ and $S_i$ is a Poisson process with rate $\mu_i$. The sample paths of these processes are assumed to be right continuous on the time interval $[0, \infty)$ and to have finite left limits on $(0, \infty)$. The process $T = (T_i : i \in \mathcal{I})$ is a continuous, $\mathbf{I}$-dimensional, nondecreasing process such that for each $i$, $T_i(t)$ represents the cumulative amount of (bandwidth) capacity allocated to route $i$ up to time $t$; it is given by

$$(89) \qquad T_i(t) = \int_0^t \Lambda_i(N(s))\, ds \qquad \text{for } i \in \mathcal{I},\ t \geq 0,$$

where the function $\Lambda(\cdot) = (\Lambda_i(\cdot) : i \in \mathcal{I})$ is defined as in Section 3 using the optimization problem (1)–(3). For each $j \in \mathcal{J}$, let

$$(90) \qquad U_j(t) = C_j t - \sum_{i \in \mathcal{I}} A_{ji} T_i(t) = C_j t - (AT(t))_j, \qquad t \geq 0.$$

The process $U_j$ is continuous and nondecreasing, and $U_j(t)$ represents the cumulative unused capacity for resource $j$ up to time $t$.

We develop a fluid model, which can be thought of as a formal law of large numbers approximation for the stochastic model. For this, we consider a sequence of stochastic systems, indexed by $r \to \infty$, in which the initial conditions may change with $r$, but $E$ and $S$ are kept fixed. The processes $N^r, U^r, T^r, E, S$, in the $r$th system are rescaled with law of large numbers scaling to yield new processes $\overline{N}^r, \overline{U}^r, \overline{T}^r, \overline{E}^r, \overline{S}^r$:

$$(91) \qquad \overline{N}^r(t) = N^r(rt)/r,$$

$$(92) \qquad \overline{U}^r(t) = U^r(rt)/r,$$

$$(93) \qquad \overline{T}^r(t) = T^r(rt)/r,$$

$$(94) \qquad \overline{E}^r(t) = E(rt)/r,$$

$$(95) \qquad \overline{S}^r(t) = S(rt)/r$$

for each $t \geq 0$.

Theسsample paths of $(\overline{T}^r, \overline{E}^r, \overline{S}^r, \overline{N}^r, \overline{U}^r)$ lie in the space $D^{4\mathbf{I}+\mathbf{J}}$ of functions defined from $[0, \infty)$ into $\mathbb{R}^{4\mathbf{I}+\mathbf{J}}$ that are right continuous on $[0, \infty)$ and have finite limits from the left on $(0, \infty)$. This space is endowed with the usual Skorokhod $J_1$-topology (cf. [10]). We say that a subsequence (which may be the whole sequence) of $\{(\overline{T}^r, \overline{E}^r, \overline{S}^r, \overline{N}^r, \overline{U}^r)\}$ is tight (resp. converges weakly) if the probability measures induced on $D^{4\mathbf{I}+\mathbf{J}}$ by the subsequence are tight (resp. converge weakly). We say that the subsequence is $\mathbf{C}$-tight if it is tight and any weak limit point has continuous paths almost surely.

To obtain the fluid model, we formally pass to the limit as $r \to \infty$ in the rescaled versions of (88)–(90). In fact, we have the following theorem. For this, recall the definition of a regular point given before Definition 5.1.

THEOREM B.1. *Suppose that $\{\overline{N}^r(0)\}$ converges in distribution as $r \to \infty$ to a random variable taking values in $\mathbb{R}_+^{\mathbf{I}}$. Then, the sequence $\{(\overline{T}^r, \overline{E}^r, \overline{S}^r, \overline{N}^r, \overline{U}^r)\}$ is $\mathbf{C}$-tight, and any weak limit point $(\overline{T}, \overline{E}, \overline{S}, \overline{N}, \overline{U})$ of this sequence (obtained by taking a weak limit as $r \to \infty$ along a suitable subsequence), almost surely satisfies the following for all $t \geq 0$: $\overline{E}(t) = \nu t$, $\overline{S}(t) = \mu t$,*

$$\overline{N}_i(t) = \overline{N}_i(0) + \nu_i t - \mu_i \overline{T}_i(t) \geq 0, \qquad i \in \mathcal{I}, \tag{96}$$

$$\overline{U}_j(t) = C_j t - (A\overline{T}(t))_j \geq 0, \qquad j \in \mathcal{J}, \tag{97}$$

*where $\overline{U}$ is nondecreasing and continuous, $\overline{T}$ is uniformly Lipschitz continuous with a Lipschitz constant bounded by $C^* = \max_{j \in \mathcal{J}} C_j$, and $\overline{T}(0) = 0$. Moreover, for a.e. sample point $\omega$, $\overline{T}(\cdot, \omega)$ is absolutely continuous and for each regular point $t$ for $\overline{T}(\cdot, \omega)$, for each $i \in \mathcal{I}$,*

$$\frac{d}{dt}\overline{T}_i(t, \omega) = \begin{cases} \Lambda_i(\overline{N}(t, \omega)), & \text{if } i \in \mathcal{I}_+(\overline{N}(t, \omega)), \\ \rho_i, & \text{if } i \in \mathcal{I}_0(\overline{N}(t, \omega)), \end{cases} \tag{98}$$

*and* (13) *and* (14) *hold at $t$ for each $i \in \mathcal{I}$ and $j \in \mathcal{J}$, with $n(\cdot) = \overline{N}(\cdot, \omega)$ there.*

PROOF. For each value of $r$ and $t \geq 0$,

$$\overline{T}^r(t) = \frac{1}{r} \int_0^{rt} \Lambda(N^r(s)) \, ds = \int_0^t \Lambda(\overline{N}^r(s)) \, ds, \tag{99}$$

where we have used the scaling property of Lemma A.2. Hence, since $\Lambda(\cdot)$ is uniformly bounded by $C^*$, for each $r$, $t \to \overline{T}^r(t)$ is uniformly Lipschitz continuous with a Lipschitz constant bounded by $C^*$. It follows immediately that the sequence of processes $\{\overline{T}^r(\cdot)\}$ is $\mathbf{C}$-tight. On combining this with the assumed convergence in distribution of $\{\overline{N}^r(0)\}$, the functional weak law of large numbers for $(E, S)$, (88) and (90), we see that $\{(\overline{T}^r, \overline{E}^r, \overline{S}^r, \overline{N}^r, \overline{U}^r)\}$ is $\mathbf{C}$-tight, and any weak limit point $(\overline{T}, \overline{E}, \overline{S}, \overline{N}, \overline{U})$ of this sequence has the following properties: almost surely for all $t \geq 0$, $\overline{E}(t) = \nu t$, $\overline{S}(t) = \mu t$ and (96) and (97) hold, where $\overline{U}$ is nondecreasing and continuous, and $\overline{T}$ is uniformly Lipschitz continuous with a Lipschitz constant bounded by $C^*$, and $\overline{T}(0) = 0$.

To prove the remaining properties claimed in the theorem, without loss of generality, using the Skorokhod representation theorem and the fact that the sample paths of $\overline{T}, \overline{E}, \overline{S}, \overline{N}, \overline{U}$ are continuous with probability one, we may assume that all stochastic processes under consideration are defined on the same probability space and that almost surely

$$\{(\overline{T}^r, \overline{E}^r, \overline{S}^r, \overline{N}^r, \overline{U}^r)\} \to (\overline{T}, \overline{E}, \overline{S}, \overline{N}, \overline{U}) \qquad \text{u.o.c.} \tag{100}$$

as $r \to \infty$ through a sequence $\{r_k\}_{k=1}^\infty$. Here u.o.c. denotes uniformly on compact time intervals.

Fix a sample point $\omega$ such that the convergence in (100) holds for $\omega$. Then, $\overline{T}(\cdot,\omega)$ is absolutely continuous and differentiable at almost every time. Fix a regular point $t$ for $\overline{T}(\cdot,\omega)$. This is the same as a regular point for $\overline{N}(\cdot,\omega)$, and $\overline{U}(\cdot,\omega)$ is differentiable there. For $i \in \mathcal{I}_0(\overline{N}(t,\omega))$, $\overline{N}_i(t,\omega) = 0$ and since $\overline{N}_i(\cdot,\omega) \geq 0$, we must have $\frac{d}{dt}\overline{N}_i(t,\omega) = 0$, and hence, by (96), $\frac{d}{dt}\overline{T}_i(t,\omega) = \rho_i = \nu_i/\mu_i$. For $i \in \mathcal{I}_+(\overline{N}(t,\omega))$, $\overline{N}_i(t,\omega) > 0$ and by the continuity of $\overline{N}(\cdot,\omega)$, there is $\varepsilon > 0$ (depending on $t$ and $\omega$) such that $\overline{N}_i(s,\omega) > 0$ for all $s \in [t, t+\varepsilon]$. By (100), $\overline{N}^{r_k}(\cdot,\omega) \to \overline{N}(\cdot,\omega)$ uniformly on $[t, t+\varepsilon]$ as $k \to \infty$, and it then follows from the continuity property of Lemma A.3 that $\Lambda_i(\overline{N}^{r_k}(\cdot,\omega)) \to \Lambda_i(\overline{N}(\cdot,\omega))$ uniformly on $[t, t+\varepsilon]$ as $k \to \infty$. Hence, by passing to the limit in (99) we obtain

$$(101) \quad \overline{T}_i(s,\omega) - \overline{T}_i(t,\omega) = \int_t^s \Lambda_i(\overline{N}(u,\omega))\,du \qquad \text{for all } s \in [t, t+\varepsilon].$$

It follows, that $\frac{d}{dt}\overline{T}_i(t,\omega) = \Lambda_i(\overline{N}(t,\omega))$, since $\Lambda_i(\cdot)$ is continuous on $\{n \in \mathbb{R}_+^{\mathbf{I}} : n_i > 0\}$, by Lemma A.3. Combining the above, we see that (98) holds at a regular point $t$ for $\overline{T}(\cdot,\omega)$. Moreover, for each $j \in \mathcal{J}$, since $\overline{U}_j(\cdot,\omega)$ is a nondecreasing function, we have

$$(102) \qquad 0 \leq \frac{d}{dt}\overline{U}_j(t,\omega) = C_j - \left(A\frac{d}{dt}\overline{T}(t,\omega)\right)_j.$$

Combining (98) and (102) with (96), we see that $n(\cdot) = \overline{N}(\cdot,\omega)$ satisfies (13) and (14) for each $i \in \mathcal{I}$ and $j \in \mathcal{J}$, at each regular point $t$ for $\overline{N}(\cdot,\omega)$ [which is the same as a regular point for $\overline{T}(\cdot,\omega)$]. $\square$

REMARK B.1. Theorem B.1 shows that there exists a fluid model solution for each possible starting point $n(0) \in \mathbb{R}_+^{\mathbf{I}}$. Note, however, that we are not assuming uniqueness of fluid model solutions given the starting state.

**Acknowledgment.** The authors thank Maury Bramson for helpful conversations concerning his related work.

## REFERENCES

[1] BEN FREDJ, S., BONALD, T., PROUTIERE, A., REGNIE, G. and ROBERTS, J. (2001). Statistical bandwidth sharing: A study of congestion at flow level. In *Proc. ACM Sigcomm 2001*.
[2] BONALD, T. and MASSOULIÉ, L. (2001). Impact of fairness on Internet performance. In *Proc. ACM Sigmetrics 2001*.
[3] BRAMSON, M. (1996). Convergence to equilibria for fluid models of FIFO queueing networks. *Queueing Syst. Theory Appl.* **22** 5–45. MR1393404

[4] BRAMSON, M. (1997). Convergence to equilibria for fluid models of head-of-the-line proportional processor sharing queueing networks. *Queueing Syst. Theory Appl.* **23** 1–26. MR1433762

[5] BRAMSON, M. (1998). State space collapse with application to heavy traffic limits for multiclass queueing networks. *Queueing Syst. Theory Appl.* **30** 89–148. MR1663763

[6] BRAMSON, M. and DAI, J. G. (2001). Heavy traffic limits for some queueing networks. *Ann. Appl. Probab.* **11** 49–90. MR1825460

[7] BU, T. and TOWSLEY, D. (2001). Fixed point approximation for TCP behavior in an AQM network. In *Proc. ACM Sigmetrics 2001*.

[8] CLARKE, F. H. (1990). *Optimization and Non Smooth Analysis*. SIAM, Philadelphia. MR1058436

[9] DE VECIANA, G., LEE, T.-J. and KONSTANTOPOULOS, T. (2001). Stability and performance analysis of networks supporting elastic services. *IEEE/ACM Transactions on Networking* **9** 2–14.

[10] ETHIER, S. N. and KURTZ, T. G. (1986). *Markov Processes*: *Characterization and Convergence*. Wiley, New York. MR838085

[11] FLOYD, S. and FALL, K. (1999). Promoting the use of end-to-end congestion control in the Internet. *IEEE/ACM Transactions on Networking* **7** 458–472.

[12] GIBBENS, R. J., SARGOOD, S. K., VAN EIJL, C., KELLY, F. P., AZMOODEH, H., MACFADYEN, R. N. and MACFADYEN, N. W. (2000). Fixed-point models for the end-to-end performance analysis of IP networks. In *13th ITC Specialist Seminar*: *IP Traffic Measurement*, *Modeling and Management*.

[13] HARRISON, J. M. (2000). Brownian models of open processing networks: Canonical representation of workload. *Ann. Appl. Probab.* **10** 75–103. [Correction (2003) **13** 390–393.] MR1765204

[14] HARRISON, J. M. (2003). A broader view of Brownian networks. *Ann. Appl. Probab.* **13** 1119–1150. MR1994047

[15] LUENBERGER, D. G. (1969). *Optimization by Vector Space Methods*. Wiley, New York. MR238472

[16] MATHIS, M., SEMKE, J., MAHDAVI, J. and OTT, T. (1997). The macroscopic behaviour of the TCP congestion avoidance algorithm. *Computer Communication Review* **27** 67–82.

[17] MO, J. and WALRAND, J. (2000). Fair end-to-end window-based congestion control. *IEEE/ACM Transactions on Networking* **8** 556–567.

[18] PADHYE, J., FIROIU, V., TOWSLEY, D. and KUROSE, J. (2000). Modeling TCP Reno performance: A simple model and its empirical validation. *IEEE/ACM Transactions on Networking* **8** 133–145.

[19] ROBERTS, J. and MASSOULIÉ, L. (2000). Bandwidth sharing and admission control for elastic traffic. *Telecomunication Systems* **15** 185–201.

[20] ROUGHAN, M., ERRAMILLI, A. and VEITCH, D. (2001). Network performance for TCP networks, Part I: Persistent Sources. In *Proc. International Teletraffic Congress* 24–28.

[21] WILLIAMS, R. J. (1998). Diffusion approximations for open multiclass queueing networks: Sufficient conditions involving state space collapse. *Queueing Syst. Theory Appl.* **30** 27–88. MR1663759

Statistical Laboratory  
Centre for Mathematical Sciences  
University of Cambridge  
Wilberforce Road  
Cambridge CB3 0WB  
United Kingdom  
e-mail: F.P.Kelly@statslab.cam.ac.uk

Department of Mathematics  
University of California  
San Diego  
La Jolla, California 92093-0112  
USA  
e-mail: williams@math.ucsd.edu